\documentclass[a4paper]{amsart}

\usepackage{amsmath,amssymb}
\usepackage[dvipsnames]{xcolor}
\usepackage{pictexwd}
\usepackage{graphicx}

\makeatletter
\def\l@subsection{\@tocline{2}{0pt}{3.5pc}{0pc}{}}
\makeatother

\newtheorem{theorem}{Theorem}[section]
\newtheorem{proposition}[theorem]{Proposition}
\newtheorem{corollary}[theorem]{Corollary}
\newtheorem{lemma}[theorem]{Lemma}
\newtheorem{algorithm}[theorem]{Algorithm}

\newtheorem{preremark}[theorem]{Remark}
\newtheorem{predefinition}[theorem]{Definition}
\newtheorem{preexample}[theorem]{Example}
\newtheorem{prenotation}[theorem]{Notation}
\newtheorem{preconjecture}[theorem]{Conjecture}

\newenvironment{remark}{\begin{preremark}\rm}{\end{preremark}}
\newenvironment{definition}{\begin{predefinition}\rm}
{\end{predefinition}}
\newenvironment{example}{\begin{preexample}\rm}{\end{preexample}}

\newenvironment{notation}{\begin{prenotation}\rm}{\end{prenotation}}

\allowdisplaybreaks
\def\OO{{\mathcal{O}}}

\def\AA{\mathbb{A}}

\def\QQ{\mathbb{Q}}

\newcommand{\M}{{\mathfrak{M}}}

\newcommand{\m}{{\mathfrak{m}}}

\let\epsilon=\varepsilon

\def\phi{{\varphi}}
\let\Psi=\varPsi
\let\Phi=\varPhi
\let\theta=\vartheta
\let\rho=\varrho

\def\LT{\mathop{\rm LT}\nolimits}

\def\ND{\mathop{\rm ND}\nolimits}
\def\AR{\mathop{\rm AR}\nolimits}
\def\GL{\mathop{\rm GL}\nolimits}

\def\Mat{\mathop{\rm Mat}\nolimits}

\def\Supp{\mathop{\rm Supp}\nolimits}
\def\Spec{\mathop{\rm Spec}\nolimits}

\def\Cot{\mathop{\rm Cot}\nolimits}

\def\GFan{\mathop{\rm GFan}\nolimits}
\def\LTGFan{\mathop{\rm LTGFan}\nolimits}

\def\Ker{\mathop{\rm Ker}\nolimits}

\newcommand{\Lin}{\mathop{\rm Lin}\nolimits}

\newcommand{\BO}{\mathbb{B}_{\mathcal{O}}}

\newcommand{\cupdot}{\mathbin{\mathaccent\cdot\cup}}
\renewcommand{\circle}{{\circ}}

\let\To=\longrightarrow

\def\tr{^{\,\rm tr}}
\def\tfrac #1#2{{\textstyle\frac{#1}{#2}}}
\def\tsum_#1^#2{{\textstyle\sum\limits_{#1}^{#2}}}

\definecolor{red}{rgb}{1.0, 0.0, 0.0}

\def\cocoa{\mbox{\rm
  C\kern-.13em o\kern-.07 em C\kern-.13em o\kern-.15em A}}
\def\apcocoa{\mbox{\rm
A\kern-0.13em p\kern -0.07em C\kern-.13em o\kern-.07 em C\kern-.13em
o\kern-.15em A}}

\begin{document}

\title{Re-embeddings of Affine Algebras Via Gr\"obner Fans of Linear Ideals}

\author{Martin Kreuzer}
\address{Fakult\"at f\"ur Informatik und Mathematik, Universit\"at
Passau, D-94030 Passau, Germany}
\email{Martin.Kreuzer@uni-passau.de}

\author{Le Ngoc Long}
\address{Fakult\"at f\"ur Informatik und Mathematik, Universit\"at
Passau, D-94030 Passau, Germany and Department of Mathematics,
University of Education - Hue University, 34 Le Loi Street, Hue City, Vietnam}
\email{lengoclong@dhsphue.edu.vn}

%    Information for third author
\author{Lorenzo Robbiano}
\address{Dipartimento di Matematica, Universit\`a di Genova,
Via Dodecaneso 35,
I-16146 Genova, Italy}
\email{lorobbiano@gmail.com}

\date{\today}

\begin{abstract}
Given an affine algebra $R=K[x_1,\dots,x_n]/I$ over a field~$K$, where~$I$
is an ideal in the polynomial ring $P=K[x_1,\dots,x_n]$, we examine the task of
effectively calculating re-embeddings of~$I$, i.e., of presentations $R=P'/I'$ such that
$P'=K[y_1,\dots,y_m]$ has fewer indeterminates. For cases when the number of indeterminates~$n$
is large and Gr\"obner basis computations are infeasible, we have introduced the
method of $Z$-separating re-embeddings in~\cite{KLR2} and~\cite{KLR3}.
This method tries to detect polynomials of a special shape in~$I$ which allow us
to eliminate the indeterminates in the tuple~$Z$ by a simple substitution process.
Here we improve this approach by showing that suitable candidate tuples~$Z$
can be found using the Gr\"obner fan of the linear part of~$I$. Then we describe
a method to compute the Gr\"obner fan of a linear ideal, and we improve this
computation in the case of binomial linear ideals using a cotangent equivalence relation.
Finally, we apply the improved technique in the case of the defining ideals of border
basis schemes.
\end{abstract}



\maketitle

%\bigskip
%\tableofcontents

%%%%%%%%%%%%%%%%%%%%%%%%%%%%%%%%%%%%%%
%
%  Section 1: Introduction
%
%%%%%%%%%%%%%%%%%%%%%%%%%%%%%%%%%%%%%%

\section{Introduction}

A finitely generated algebra~$R$ over a field~$K$ is also called an {\it affine $K$-algebra}.
In order to analyse an affine $K$-algebra and to perform basic operations effectively, we usually
assume that the algebra is given by generators and relations, i.e., that $R=P/I$
where $P=K[x_1,\dots,x_n]$ is a polynomial ring over~$K$ and~$I$ is an explicitly given ideal in~$P$.
The question of whether we can actually perform the intended calculations, many of which are based
on the theory of Gr\"obner bases, depends chiefly on the number~$n$ of indeterminates involved
in the presentation $R=P/I$. If the number~$n$ is too large, essentially all computations, with the
possible exception of those which require solely linear algebra, become prohibitively expensive. 
Thus it is an important task to find better presentations $R \cong P'/I'$ where $P'=K[y_1,\dots,y_m]$ 
is a polynomial ring in fewer indeterminates and~$I'$ is an ideal in~$P'$. 

In the language of Algebraic Geometry, we are given an affine scheme $\Spec(R)$ embedded into
an affine space $\mathbb{A}^n_K$, and we are looking for a re-embedding of it into a lower
dimensional affine space $\mathbb{A}^m_K$. In this setting, the topic has a long history, beginning
with the classical result that a smooth variety of dimension~$d$ can be embedded into $\mathbb{A}^{2d+1}$
if the field~$K$ is infinite. In the fundamental paper~\cite{Sri}, this bound is generalized to
a bound for arbitrary affine schemes over infinite fields. The usual way to achieve the desired
improvements of an embedding in algebraic geometry is by using {\it generic projections}.
The obstruction is then mainly given by the secant variety of the scheme (see~\cite[Section~7]{Hol}).
If we intend to perform actual computer calculations, computing projections corresponds to
calculating elimination ideals, possibly after a change of coordinates. Herein lies the heart of the 
problem: finding suitable projections and calculating the resulting presentations $R \cong P'/I'$
is traditionally done by computing Gr\"obner bases with respect to elimination term orderings, and this
is usually infeasible when the number of indeterminates is large.

In two previous papers~\cite{KLR2} and~\cite{KLR3}, the authors introduced and developed
a new method for re-embedding affine algebras based on tuples of {\it separating indeterminates}.
This method can be applied to any finitely generated $K$-algebra $R=P/I$ if~$I$ is an ideal 
contained in the maximal ideal $\M = \langle x_1,\dots,x_n \rangle$ of~$P$. 
If this hypothesis is not satisfied from the outset, it may be necessary to know a $K$-rational
point of $\Spec(R)$ and to perform a suitable change of coordinates.
Then, for a tuple $Z=(z_1,\dots,z_s)$ of indeterminates
in~$P$, we say that~$I$ is {\it $Z$-separating} if there exist a tuple of non-zero 
polynomials $(f_1,\dots,f_s)$ of~$I$ and a term ordering~$\sigma$ such that
$\LT_\sigma(f_i)=z_i$ for $i=1,\dots,s$. Given a tuple~$Z$ such that~$I$ is $Z$-separating, 
we can eliminate the indeterminates in~$Z$ and get a {\it $Z$-separating
re-embedding} $\Phi:\; P/I \longrightarrow K[Y] / (I\cap K[Y])$ where $Y= \{x_1,\dots,x_n\}
\setminus \{z_1,\dots,z_s\}$. Here the actual elimination can be carried out by 
interreducing $(f_1,\dots,f_s)$ such that they become {\it coherently separating} and
performing substitutions in the remaining generators of~$I$. Of course, in this way we will usually
end up with a system of generators of the elimination ideal, but not a Gr\"obner basis with respect
to any term ordering. 

Notice that the question of whether an ideal~$I$ in~$P$ is $Z$-separating for some tuple of
indeterminates~$Z$ is related to
the famous {\it Epimorphism Problem} in Affine Geometry (see for instance~\cite{vdE}, 
\cite{Dre}, or~\cite{Gup}). However, we are looking for even more special generators of~$I$
than required in this problem, and therefore we generally cannot expect~$I$ to be $Z$-separating 
for some tuple~$Z$. Our main focus is to try to find such tuples~$Z$ computationally, even if we 
cannot know whether they exist for a given ideal~$I$, because for larger numbers of indeterminates 
this is our best chance to effectively calculate good, and possibly optimal, re-embeddings.

Using the method of $Z$-separating re-embeddings, the task of finding 
good re-embeddings is split into two steps:

\noindent (1) Find a suitable candidate tuple of indeterminates~$Z$. 

\noindent (2) Check if this tuple~$Z$ really works.

For the second step we can use a method based on Linear Programming Feasibility (LPF),
if we have a candidate tuple $(f_1,\dots,f_s)$ of polynomials $f_i\in I$ available (see~\cite[Section~3]{KLR3}).
Other algorithms for this task, which are efficient but not guaranteed to succeed,
are in development. The current paper is chiefly concerned with the
first step of finding suitable candidate tuples~$Z$. A first idea of how to restrict
the number of tuples one has to consider was introduced in~\cite[Section~5]{KLR3},
namely the idea is to use the {\it $Z$-restricted Gr\"obner fan} of~$I$. However, this approach may still be
computationally heavy, and this is where the first new result of the current paper comes in:
it turns out that it is enough to consider the Gr\"obner fan of the {\it linear part}
$\Lin_\M(I)$ of~$I$ (see Propositions~\ref{prop-cohSepPoly} and~\ref{prop-charLinZ}). 
The linear part of~$I$ is the linear ideal generated by the linear parts
of an arbitrary system of generators of~$I$, and thus easy to compute. 

In Section~2 we begin by recalling and extending the basic theory of $Z$-separating re-embeddings. 
For the later applications, it is necessary to pay close attention
to the questions of when the newly found re-embeddings are {\it optimal}, i.e., 
when the number~$m$ of indeterminates of~$P'$ is the minimal possible one (see Corollary~\ref{cor-cohSepPoly}), 
and when the re-embedded ring is actually a polynomial ring, i.e., when $I'=\langle 0\rangle$
and $\Spec(R)$ is isomorphic to an affine space (see Corollary~\ref{cor-affinescheme}).

Thus we are led to examine the task of calculating the Gr\"obner fan of a linear ideal $I_L$ in
Section~3. First we reduce the problem to check whether a given
tuple of indeterminates~$Z$ consists of the leading terms of a reduced Gr\"obner basis
of~$I_L$ to a linear algebra computation (see Proposition~\ref{prop-charLinZ}). Then
we reduce the calculation of the Gr\"obner fan of~$I_L$
to the task of finding the maximal minors of a matrix (see Theorem~\ref{thm-bijection}). 
We note that this step can also be tackled via any 
method for computing the bases of a linear matroid (see Remark~\ref{rem-compGFan}.b).
Combining all steps, we present an algorithm for finding $Z$-separating re-embeddings
using the Gr\"obner fan of $\Lin_\M(I)$ in Section~4. Notice that Algorithm~\ref{alg-Emb} also allows us
to find optimal $Z$-separating re-embeddings, if they exist.

As the calculation of the maximal minors of a matrix required by Theorem~\ref{thm-bijection} 
could still be quite demanding, we look at a practically relevant special case in Section~5,
namely the case of binomial linear ideals. In this case we define an equivalence relation
on the indeterminates in $X=(x_1,\dots,x_n)$, called {\it cotangent equivalence} and based
on the equality of $K\cdot \bar{x}_i$ and $K\cdot \bar{x}_j$ in the cotangent
space $\Cot_\m(R)= \m/\m^2$, where $\m = \langle \bar{x}_1,\dots,\bar{x}_n\rangle$ is the 
maximal ideal of~$R$ generated by the residue classes of the indeterminates.
Using this equivalence relation, the indeterminates in~$X$ can be divided into three types:
{\it basic} indeterminates which can never be a part of any separating tuple~$Z$,
{\it trivial} indeterminates which can be part of such a tuple~$Z$, and {\it proper} 
equivalence classes of indeterminates for which a proper subset can be put into~$Z$.
With the help of this equivalence relation we describe the leading term Gr\"obner fan
of a linear ideal explicitly (see Theorem~\ref{thm-shapeofSsigma}), apply it to the
classification of possible separating tuples~$Z$ (see Theorem~\ref{thm-classifycij}), and
spell out an explicit efficient algorithm to compute those tuples (see Algorithm~\ref{alg-compEmb}).

In the final section we apply the new theory to the case which prompted its development
in the first place, namely to the coordinate rings of border basis schemes. 
These schemes are important moduli spaces in Algebraic Geometry which parametrize 0-dimensional
polynomial ideals. Their vanishing ideals have exactly the correct structure to make
$Z$-separating re-embeddings work: they are contained in the maximal ideal generated by the indeterminates,
their linear part is a binomial linear ideal, and the cotangent equivalence classes can be
calculated quickly (using the algorithm given in~\cite{KSL}). 
In this setting the indeterminates of the underlying polynomial ring can be classified
further into {\it interior} indeterminates and {\it rim} indeterminates, and we are
able to provide additional information about their distribution in the cotangent equivalence classes
(see Theorem~\ref{thm-RimIndets}). Our final Example~\ref{ex-Section6} shows the
machinery of $Z$-separating re-embeddings at work and verifies some 
claims in~\cite[Remark~7.5.3]{Hui1}.

All algorithms mentioned in this paper were implemented in the computer 
algebra system CoCoA (see~\cite{CoCoA}) and are available as a package
for ApCoCoA (see~\cite{ApCoCoA}) on the first author's web 
page\footnote{\tt https://symbcomp.fim.uni-passau.de/en/symbolic-computation/projects}
This package can be applied to perform the calculations underlying most examples. 
Their use was essential in the discovery of properties and features which eventually 
evolved into theorems or disproved previous conjectures. 
The general notation and definitions in this paper follow~\cite{KR1} and~\cite{KR2}.

\bigskip\bigbreak
%%%%%%%%%%%%%%%%%%%%%%%%%%%%%%%%%%%%%%%%%%%%%%%%%%%%%%%%%%%%%%%%%%%
%
% Section 2: Z-Separating Re-embeddings
%
%%%%%%%%%%%%%%%%%%%%%%%%%%%%%%%%%%%%%%%%%%%%%%%%%%%%%%%%%%%%%%%%%%%

\section{Z-Separating Re-embeddings}
\label{Z-Separating Re-embeddings}

In this paper we let $K$ be an arbitrary field, let $P=K[x_1,\dots,x_n]$,
and let $\M=\langle x_1, \dots, x_n \rangle$. The tuple formed by the indeterminates 
of~$P$ is denoted by $X =(x_1, \dots, x_n)$. Moreover, let $1\le s\le n$, 
let $z_1,\dots, z_s$ be pairwise distinct indeterminates in~$X$, 
and let $Z =(z_1, \dots, z_s)$.
Denote the remaining indeterminates by 
$\{y_1,\dots,y_{n-s}\} = \{x_1,\dots,x_n\} \setminus \{z_1,\dots,z_s\}$,
and let $Y =(y_1,\dots, y_{n-s})$. Committing a slight abuse of notation, 
we shall also write $Y = X \setminus Z$.
The monoid of terms in~$P$ is denoted by $\mathbb{T}^n = \{ x_1^{\alpha_1} \cdots
x_n^{\alpha_n} \mid \alpha_i \ge 0\}$. 
Given a term ordering $\sigma$ on~$\mathbb{T}^n$,
its restriction to $\mathbb{T}(y_1, \dots, y_{n-s})$ is denoted by~$\sigma_Y$.

Recall that an algebra of type $R=P/I$ where $P=K[x_1,\dots,x_n]$ is a polynomial ring over a field $K$ 
and~$I$ is a proper ideal in~$P$ is called an {\bf affine $K$-algebra}. In this setting,
re-embeddings are defined as follows.

\begin{definition}\label{def:re-embedding}
Let $P=K[x_1,\dots,x_n]$ be a polynomial ring over a field~$K$, let~$I$ 
be a proper ideal in~$P$, and let $R=P/I$.
\begin{enumerate}
\item[(a)] A $K$-algebra isomorphism $\Psi:\; R \longrightarrow P'/I'$,
where $P'$ is a polynomial ring over~$K$ and~$I'$ is an ideal in~$P'$, is
called a {\bf re-embedding} of~$I$.

\item[(b)] A re-embedding $\Psi:\; R \longrightarrow P'/I'$ of~$I$ is called
{\bf optimal} if every \hbox{$K$-al}\-ge\-bra isomorphism $R \longrightarrow P''/I''$ 
with a polynomial ring~$P''$ over~$K$ and an ideal~$I''$ in~$P''$ 
satisfies the inequality $\dim(P'')\ge \dim(P')$.
\end{enumerate}
\end{definition}

In~\cite{KLR2} and~\cite{KLR3}, the authors examined re-embeddings of affine $K$-algebras,
i.e., isomorphisms with presentations requiring fewer $K$-algebra generators. 
For these techniques to work, we need to assume that the given ideal~$I$ is contained 
in a linear maximal ideal of~$P$. As explained in \cite[Section~1]{KLR2}, we can then 
perform a linear change of coordinates and assume that~$I$ is contained in
$\M= \langle x_1,\dots,x_n\rangle$.
In particular, it was shown that the following situation leads to such re-embeddings.

\begin{definition}\label{def-Zseparating}
Let $I$ be an ideal in~$P$ with $I\subseteq\M$, and let $Z=(z_1,\dots,z_s)$
be a tuple of distinct indeterminates in~$X$. 
\begin{enumerate}
\item[(a)]
We say that the ideal~$I$ is {\bf $Z$-separating} if there exist a term ordering~$\sigma$ on~$\mathbb{T}^n$ 
and $f_1,\dots,f_s \in I \setminus \{0\}$ such that
$\LT_\sigma(f_i)= z_i$ for $i=1,\dots,s$. In this situation~$\sigma$ is called a 
{\bf $Z$-separating term ordering} for~$I$, and the tuple $(f_1,\dots,f_s )$ is called
a {\bf $Z$-separating tuple}.

\item[(b)] The ideal~$I$ is called {\bf coherently $Z$-separating} if it contains a
$Z$-sepa\-rat\-ing tuple $(f_1,\dots,f_s)$ such that for $i\ne j$ the indeterminate~$z_i$
does not divide any term in the support of~$f_j$.
\end{enumerate}
\end{definition}

Given a $Z$-separating term ordering~$\sigma$ for~$I$,
the reduced $\sigma$-Gr\"obner basis of~$I$ is of the form 
$G = \{ z_1 - h_1,\allowbreak \dots,\, z_s - h_s, g_1,\, \dots,\, g_r\}$
with $h_i, g_j \in K[Y]$. 
In this case the $K$-algebra homomorphism $\Phi: P/I \longrightarrow K[Y] / 
(I \cap K[Y])$ given by $\Phi(\bar{x}_i) = \bar{x}_i$ for 
$x_i\in Y$ and $\Phi(\bar{x}_i) = \bar{h}_j$
for $x_i = z_j \in Z$ is an isomorphism of $K$-algebras. It is called the
{\bf $Z$-separating re-embedding} of~$I$ (see~\cite[Theorem~2.13]{KLR2}).
Notice that the map~$\Phi$ is a re-embedding of~$I$ such that the new polynomial ring
$K[Y]$ involves fewer indeterminates, and the size of~$Z$ measures the improvement
$\#Z = \#X - \#Y$ we achieved. Geometrically, the original variety can be viewed
as the graph of the functions $h_1,\dots,h_s$ over the re-embedded variety.

For the choice of a $Z$-separating term ordering, we have the following observation.

\begin{proposition}\label{prop-sepZ}
Let~$I$ be an ideal in~$P$ which is contained in~$\M$, and
let $Z=(z_1,\dots,z_s)$ be a tuple of distinct indeterminates in~$X$.
Then the following conditions are equivalent.
\begin{enumerate}
\item[(a)] The ideal~$I$ is $Z$-separating.

\item[(b)] For every elimination ordering~$\sigma$ for~$Z$,
we have $\langle Z\rangle \subseteq \LT_\sigma(I)$.

\item[(c)] There exists an elimination ordering~$\sigma$ for~$Z$ such that
$\langle Z\rangle \subseteq \LT_\sigma(I)$.
\end{enumerate}
\end{proposition}

\begin{proof}
To show that~(a) implies~(b), we note that, by definition, there exists
a $Z$-separating term ordering~$\sigma$ for~$I$. By~\cite[Remark~4.3]{KLR3}, it
follows that any elimination ordering for~$Z$ is then also a 
$Z$-separating term ordering for~$I$. Now the claim follows 
from~\cite[Proposition~4.2]{KLR3}.

Condition (b) obviously implies~(c), and the remaining implication
(c)$\Rightarrow$(a) follows from the definition.
\end{proof}

The following example illustrates this proposition.

\begin{example}\label{ex-GBnecessary}
Consider the ring $P=\QQ[x,y,z]$, the tuple $Z=(x)$, and the ideal 
$I = \langle f_1,\dots,f_{10}\rangle$, where 
$$
\begin{array}{lcl}
f_1 &=&  xy^2 +\frac{1}{2}y^3 -\frac{1}{2}y^2z -x^2 -\frac{1}{2}xy -y^2 +\frac{1}{2}xz +x, \cr
f_2 &=&  y^2z^2 +3y^3 -4y^2z -xz^2 -3xy +4xz  ,\cr
f_3 &=&  y^3z -xyz -y^2z +xz   ,\cr
f_4 &=&  y^4 -xy^2 -y^3 +xy   ,\cr
f_5 &=&  x^2y^2 -x^3  ,\cr
f_6 &=&  x^3 +\frac{1}{2}x^2y +xy^2 +\frac{1}{2}y^3 -\frac{1}{2}x^2z -\frac{1}{2}y^2z -x^2 -y^2 ,\cr
f_7 &=&  x^2z^2 +y^2z^2 +3x^2y +3y^3 -4x^2z -4y^2z   ,\cr
f_8 &=&  x^2yz +y^3z -x^2z -y^2z  ,\cr
f_9 &=&  x^2y^2 +y^4 -x^2y -y^3   ,\cr
f_{10} &=&  x^4 +x^2y^2.   
\end{array}
$$
At first glance, the ideal~$I$ does not appear to be $Z$-separating, even if we use
linear combinations of the generators. However, for any elimination ordering~$\sigma$ for~$Z$, the
reduced $\sigma$-Gr\"obner basis of the ideal~$I$ is $\{ x-y^2,\, y^4+y^2 \}$, 
and this proves that~$I$ is indeed $Z$-separating.
\end{example}

The above proposition provides one way to solve the following task.

\begin{remark}{\bf (Checking $Z$-Separating Tuples)}\label{rem-checkZ}\\
Given a tuple of indeterminates~$Z$, there are several methods for checking
whether the ideal~$I$ is $Z$-separating.
\begin{enumerate}
\item[(a)] Condition~(b) of the last proposition says that
we can check $\langle Z\rangle \subseteq \LT_\sigma(I)$ for any
elimination ordering~$\sigma$ for~$Z$. However, the required Gr\"obner basis computation
may be too costly, in particular if the given ideal~$I$ is not $Z$-separating.

\item[(b)] If we have a tuple of polynomials $(f_1,\dots,f_s)$ with $f_i\in I$
and want to check whether it is $Z$-separating,
we can use the methods explained in~\cite[Section~4]{KLR3}. They use 
Linear Programming Feasibility solvers and are usually very fast.

\end{enumerate}
\end{remark}

In the following, our main focus is the possibility to weed out many candidate tuples~$Z$
beforehand. In~\cite{KLR3} it was suggested to use the Gr\"obner fan of~$I$ for this purpose.
Actually, as we shall see later, it suffices to use the Gr\"obner fan of the linear part of~$I$
which we introduce now.
Recall that $\Lin_{\M}(f)$ denotes the homogeneous component of 
standard degree 1 of a polynomial $f \in \M$. It is called the {\bf linear part}
of~$f$. Given an ideal~$I$ contained in~$\M$, the $K$-vector space
$\Lin_{\M}(I) = \langle \Lin_{\M}(f) \mid f\in I \rangle_K$ 
is called the {\bf linear part} of~$I$.

In \cite[Proposition~1.6]{KLR2}, we showed
that $\Lin_{\M}(I) $ is easy to compute, since it is equal 
to $\langle \Lin_{\M}(f_1),...,\Lin_{\M}(f_s)\rangle_K$, where $\{f_1, \dots, f_s\}$ 
is any set of  generators of~$I$.
In this setting we make the following useful observation.

\begin{proposition}\label{prop-cohSepPoly}
Let $I$ be an ideal in~$P$ which is contained in~$\M$. Suppose that~$I$ is
$Z$-separating for some tuple of indeterminates $Z$ in~$X$, let~$\sigma$ 
be a $Z$-separating term ordering for~$I$, and let $Y=X \setminus Z$.
\begin{enumerate}
\item[(a)] We have $\langle Z\rangle \subseteq  \LT_\sigma(\langle \Lin_{\M}(I) \rangle)$,

\item[(b)] Let $S_\sigma$ be the set of indeterminates which generate
$\LT_\sigma(\langle \Lin_{\M}(I)\rangle )$. Then we have $Y \supseteq X \setminus S_\sigma$.  

\item[(c)]We have $Z \subseteq \LT_\tau(\langle\Lin_\M(I)\rangle)$ for every elimination ordering~$\tau$ for~$Z$.

\end{enumerate}
\end{proposition}

\begin{proof} First we show~(a). By assumption, there exists a tuple $(f_1, \dots, f_s)$ of 
polynomials in~$I$ which is  the reduced $\sigma$-Gr\"obner 
basis of $\langle f_1,\dots,f_s\rangle$ and such that $z_i = \LT_\sigma(f_i)$ 
for $i=1,\dots, s$. Then we have $z_i = \LT_\sigma(\Lin_\M(f_i))$ for $i=1,\dots, s$, 
and thus $\langle Z\rangle \subseteq \LT_\sigma(\langle\Lin_{\M}(I)\rangle)$. 

Claim~(b) follows from~(a) and the definition of~$Y$, and claim~(c) is a consequence of~(a) 
and Proposition~\ref{prop-sepZ}.
\end{proof}

Let us denote the image of~$\M$ in~$P/I$ by~$\m$. Recall that the $K$-vector space
$\Cot_{\m}(R) = \m / \m^2$ is called the {\bf cotangent space} of~$P/I$ at the origin. 

\begin{remark}\label{rem-YgeneratesCot}
As shown in~\cite[Proposition~1.8.b]{KLR2}, the canonical map $P_1 \longrightarrow \m/\m^2$
induces an isomorphism of $K$-vector spaces $P_1 / \Lin_{\M}(I) \cong \m/\m^2$. 
In the setting of the proposition, the set of the residue classes 
of the elements in $X \setminus S_\sigma$ is a $K$-basis of $P_1/\Lin_{\M}(I)$.
Therefore the residue classes of the indeterminates in $Y = X \setminus S_\sigma$
generate the cotangent space $\Cot_{\m}(R) \cong \m / \m^2$ of $P/I$ at the origin. 
\end{remark}

When we are only looking for $Z$-separating re-embeddings of~$I$ which are optimal,
the above proposition yields the following characterization.

\begin{corollary}\label{cor-cohSepPoly}
In the setting of the proposition, assume that ${s=\#Z}$ is equal to the $K$-vector space 
dimension of $\Lin_{\M}(I)$. Then the following claims hold.
\begin{enumerate}
\item[(a)] The map $\Phi:\; P/I \longrightarrow K[Y]/ (I \cap K[Y])$ is an optimal re-embedding of~$I$.

\item[(b)] We have $\langle Z\rangle = \LT_\sigma(\langle \Lin_{\M}(I) \rangle)$,
and hence $Y = X \setminus S_\sigma$. 

\item[(c)] The residue classes of the indeterminates in $Y = X \setminus Z$
form a $K$-vector space basis of the cotangent space 
$\Cot_{\m}(R) \cong \m / \m^2$ of $P/I$ at the origin. 

\end{enumerate}
\end{corollary}

\begin{proof}
Claim~(a) is a consequence of~\cite[Corollary~4.2]{KLR2}. Let us prove~(b). 
By claim~(a) of the proposition, we have 
$\langle Z\rangle \subseteq  \LT_\sigma(\langle \Lin_{\M}(I) \rangle)$.
Since we have the equality $s = \dim_K(\langle Z\rangle_K) = \dim_K(\Lin_{\M}(I))$, we deduce 
that~$Z$ minimally generates $\LT_\sigma(\langle \Lin_{\M}(I) \rangle)$.

Finally, let us prove (c). By~(b), we have $Z = S_\sigma$. As mentioned in
the preceding remark, the set of the residue classes 
of the elements in $X \setminus S_\sigma$ is a $K$-basis 
of $P_1/\Lin_{\M}(I)$. This implies the claim.
\end{proof}

An important situation in which we obtain an optimal re-embedding is described 
by the following corollary.

\begin{corollary}\label{cor-affinescheme}
In the setting of the preceding corollary, the following claims hold.
\begin{enumerate}
\item[(a)] The localization $P_{\M}/ I_{\M}$ is a regular local ring
if and only if $I \cap K[Y] = \{0\}$.
\end{enumerate}

Now assume that these conditions are satisfied.
\begin{enumerate}
\item[(b)] The map $\Phi:\; P/I \longrightarrow K[Y]$ is an isomorphism
with a polynomial ring and the scheme $\Spec(P/I)$ is isomorphic to an affine space.

\item[(c)] The reduced $\sigma$-Gr\"obner basis of~$I$ is of the form
$G = \{z_1-h_1,\dots, z_s-h_s\}$ with $h_i\in K[Y]$ and it is a minimal set of generators of~$I$.

\item[(d)] Every $Z$-separating tuple $F=(f_1,\dots,f_s)$ is a Gr\"obner basis of~$I$ 
with respect to any elimination ordering for~$Z$. It is a minimal system of generators of~$I$
and a permutable regular sequence in~$P$.

\item[(e)] The reduced Gr\"obner bases of~$I$ for all elimination orderings for~$Z$ coincide.

\end{enumerate}
\end{corollary}

\begin{proof} 
Claim (a) follows from~\cite[Proposition~6.7]{KLR3}. Parts~(b) and~(c) are
immediate consequences of~(a).

Now we show~(d). Since~$F$ is separating, we have $\langle Z\rangle = 
\LT_\sigma(\langle F_Z\rangle)$ for every elimination ordering~$\sigma$
for~$Z$. The theorem on the computation of elimination modules (cf.~\cite[Theorem~3.4.5]{KR1}) 
implies that every $\sigma$-Gr\"obner basis of~$I$ consists
of polynomials with leading terms in $\langle Z\rangle$ and polynomials in~$K[Y]$
generating $I \cap K[Y]$. By~(a), we have $I\cap K[Y] = \{0\}$, whence
it follows that~$F$ is in fact a minimal $\sigma$-Gr\"obner basis of~$I$. 
Consequently, $F$ is a system of generators of~$I$. 
Moreover, since the leading terms of $f_1,\dots,f_s$ 
form a regular sequence, also $F = (f_1,\dots,f_s)$ is a regular sequence.

Finally, to prove (e), we use the definitions and results of~\cite[Section~5]{KLR3}.
By Thm.~5.5, the map $\Gamma_Z:\; \GFan_Z(I)\longrightarrow \GFan(I\cap K[Y])$ is bijective. 
Since we have $I\cap K[Y] =\{0\}$, we obtain $\GFan(I\cap K[Y]) =\{\emptyset\}$. Therefore $\GFan_Z(I)$, 
which is not empty, has cardinality one.
\end{proof}

Let us apply this corollary in a concrete case.

\begin{example}\label{ex-mingensaffine}
Let $P = \QQ[x,y,z,w]$, let $f_1=  w^2 +x -y +3z$, $f_2 =  zw^2 +w^3 +y$, 
$f_3=w^3 -xz +yz -3z^2 +y$, and let $I =\langle f_1, f_2, f_3\rangle$. 
By substituting $y$ with $-f_2+y = -zw^2 -w^3$ in~$f_1$,
we get $f_1' = zw^2 +w^3 +w^2 +x +3z$. Then, by substituting~$y$ with  
$-zw^2 -w^3$ and~$x$ with $-f_1' +x = -zw^2 -w^3 -w^2 -3z$ in~$f_3$ 
we get~$0$. 

Consequently, we have $I = \langle f_1', f_2\rangle$ and
$(f_1', f_2)$ is clearly the reduced Gr\"obner basis of~$I$ with respect 
to every elimination ordering for $(x, y)$. 
According to the above corollary, we obtain an isomorphism $P/I\cong \QQ[z, w]$,
and $(f_1', f_2)$ is a permutable regular sequence.
\end{example}

The final example in this section shows that not all optimal embeddings fall into the area
of application of the above corollary.

\begin{example}\label{ex-isotoK[x]} 
In $\QQ[x,y]$, consider the polynomial $F = 2x^8 +8x^6y +12x^4y^2 +8x^2y^3 +2y^4 +x$ 
Then the $\QQ$-algebra homomorphism 
$\alpha: \QQ[x, y] \To \QQ[x]$ defined by $\alpha(x) = -2x^4$and $\alpha(y) = x - 4x^8$
satisfies  $\alpha(x^2 +y) = x$ as well as $\Ker(\alpha) = \langle F\rangle$.
It follows that $\bar{\alpha}: \QQ[x, y]/ \langle F\rangle \To \QQ[x]$ 
is an optimal re-embedding of $\langle F\rangle$, although~$F$ is not $x$-separating.
Unlike the case of Example~\ref{ex-GBnecessary}, there is no separating term ordering
for the ideal $\langle F\rangle$ here at all.
\end{example}

Another example of this type is the famous Koras-Russel cubic threefold
whose coordinate ring is $R=K[x,y,z,t]/\langle x + x^2y + z^2 + t^3\rangle$.
Using completely different techniques, it was shown that~$R$ is not isomorphic
to a polynomial ring in three indeterminates (see~\cite{ML} and~\cite{Cra}).

\bigskip\bigbreak
%%%%%%%%%%%%%%%%%%%%%%%%%%%%%%%%%%%%%%%%%%%%%%%%%%%%%%%%%%%
%
%  Section 3: Gr\"obner Fans of Linear Ideals
%
%%%%%%%%%%%%%%%%%%%%%%%%%%%%%%%%%%%%%%%%%%%%%%%%%%%%%%%%%%%

\section{Gr\"obner Fans of Linear Ideals}
\label{Groebner Fans of Linear Ideals}

As before, let $K$ be a field, let $P=K[x_1,\dots,x_n]$, and let~$I$ be an
ideal of~$P$ which is contained in $\M = \langle x_1,\dots,x_n \rangle$.
In the preceding section we saw that for the existence of a $Z$-separating
re-embedding of~$I$ it is necessary that we have
$\langle Z\rangle \subseteq \LT_\sigma(\langle\Lin_{\M}(I)\rangle)$
for some term ordering~$\sigma$. Here $\langle \Lin_{\M}(I)\rangle$ is a {\bf linear ideal},
i.e., an ideal generated by linear polynomials. The possible ideals
$\LT_\sigma(\langle\Lin_{\M}(I)\rangle)$ are classified by the {\bf Gr\"obner fan} (see~\cite{MR})
of $\langle\Lin_{\M}(I)\rangle$. Therefore we study Gr\"obner fans of linear ideals
in this section with a special emphasis on the task of computing them
efficiently.

In the following we let $L=(\ell_1, \dots, \ell_r)$ be a tuple of
linear forms in~$P$, and we let $I_L = \langle L\rangle$ be the linear ideal
generated by~$L$. For $i=1,\dots,r$, we write $\ell_i = a_{i1} x_1 + \cdots + a_{in} x_n$
with $a_{ij} \in K$. Then the matrix $A=(a_{ij}) \in \Mat_{r,n}(K)$ is 
called the {\bf coefficient matrix} of~$L$. In view of Proposition~\ref{prop-cohSepPoly}, 
we are interested in the condition $\langle Z\rangle \subseteq \LT_\sigma(I_L)$. It can be
rephrased as follows.

\begin{proposition}\label{prop-charLinZ}
Let $L=(\ell_1, \dots, \ell_r)$ be a tuple of $K$-linearly independent linear forms in~$P$,
let $I_L = \langle  L\rangle$, and let $A=(a_{ij})$ be the  coefficient matrix of~$L$.
Moreover, let $s\le r$, let $Z=(z_1,\dots,z_s)$ be a tuple of distinct indeterminates in $X=(x_1,\dots,x_n)$, 
and let $Y = X\setminus Z$. Then the following conditions are equivalent.
\begin{enumerate}
\item[(a)] There exists a term ordering~$\sigma$ such that 
$\langle Z\rangle \subseteq \LT_\sigma(I_L)$.

\item[(b)] The residue classes of the elements of~$Y$ generate the
$K$-vector space $P_1 / \langle L\rangle_K$.

\item[(c)] Let $i_1,\dots,i_s\in \{1,\dots,n\}$ be the indices such that
$z_j =  x_{i_j}$ for $j=1,\dots,s$. Then the columns $i_1,\dots,i_s$ of~$A$
are linearly independent.
\end{enumerate}
\end{proposition}

\begin{proof}
To show (a)$\Rightarrow$(b), we first note that $\langle Z \rangle 
\subseteq\LT_\sigma(I_L)$ implies that the canonical map $K[Y] \cong P/\langle Z\rangle
\longrightarrow P/\LT_\sigma(I_L)$ is surjective. 
Hence the residue classes of the elements of~$Y$ generate the $K$-algebra $P/\LT_\sigma(I_L)$.
By Macaulay's Basis Theorem (see~\cite[Theorem~1.5.7]{KR1}), it follows that the residue classes
of the elements of~$Y$ generate the $K$-algebra $P/I_L$.
We observe that that~$I_L$ is generated by linear forms, and therefore $P/I_L$ is isomorphic to a polynomial ring.
Thus the residue classes of a tuple of indeterminates~$Y$ are a $K$-algebra
system of generators of the ring $P/I_L$ if and only if they are a system of generators
of the $K$-vector space given by its homogeneous component $P_1/(I_L)_1$ of degree one, 
where $(I_L)_1 = \langle L\rangle_K$.

The assumption in~(b) implies that the indeterminates in~$Z$ can be expressed as linear
combinations of the indeterminates in~$Y$. Consequently, any elimination ordering~$\sigma$
for~$Z$ satisfies $\langle Z\rangle \subseteq \LT_\sigma(I_L)$. This proves (b)$\Rightarrow$(a).

Finally, we show that~(b) and~(c) are equivalent. Notice that both conditions
imply $s\le r$. The tuple $Y = X \setminus Z$
is a system of generators of the vector space $P_1 / \langle L \rangle_K$ 
if and only if $\{\ell_1,\dots,\ell_r\}$
together with~$Y$ is a system of generators of~$P_1$. This means that, if we extend~$A$ with
$n-s$ rows that are unit vectors having their non-zero entries at the positions
of the indeterminates in~$Y$, the resulting matrix~$\overline{A}$ of size $(r+n-s)\times n$ 
has the maximal rank~$n$. Now we consider the matrix $A_{(i_1,\dots,i_s)}$ 
consisting of columns $i_1,\dots,i_s$ of~$A$. By renumbering the 
indeterminates, we may assume that $i_1=1$, $\dots$, $i_s=s$, and hence
that the extended matrix is upper block triangular of the form
$$
\overline{A} \;=\; \begin{pmatrix}
A_{(i_1,\dots,i_s)} & \ast \\
0 & I_{n-s}
\end{pmatrix}
$$
where $I_{n-s}$ is the identity matrix of size $n-s$.
Now it is clear that the rows of~$\overline{A}$ generate $K^n$ if and
only if the rows of $A_{(i_1,\dots,i_s)}$ generate $K^s$, and this is equivalent
to $A_{(i_1,\dots,i_s)}$ having maximal rank~$s$. This concludes the proof
of the proposition.
\end{proof}

As a special case of the proposition, we get the following characterization of tuples~$Z$
which are leading term tuples of a marked reduced Gr\"obner basis in $\GFan(I_L)$.
Recall that a {\bf marked Gr\"obner basis} of an ideal~$J$ is a set of pairs
$$
\overline{G} \;=\; \{\, (\LT_\sigma(g_1),\, g_1), \dots, (\LT_\sigma(g_k),\,g_k)\,\}
$$
where~$\sigma$ is a term ordering and $G=\{g_1,\dots,g_k\}$ is the reduced $\sigma$-Gr\"obner
basis of~$J$. The set of all marked reduced Gr\"obner bases of~$J$ is the
{\bf Gr\"obner fan} $\GFan(J)$ of~$J$.

\begin{corollary}\label{cor-CharLinZ}
In the setting of the proposition, assume that $s=r$. Then the following
conditions are equivalent.
\begin{enumerate}
\item[(a)] There exists a term ordering~$\sigma$ such that 
$\langle Z\rangle = \LT_\sigma(I_L)$.

\item[(b)] The residue classes of the elements of~$Y$ are a $K$-basis 
of $P_1 / \langle L\rangle_K$.

\item[(c)] Let $i_1,\dots,i_s\in \{1,\dots,n\}$ be the indices such that
$z_j =  x_{i_j}$ for $j=1,\dots,s$. Then the columns $i_1,\dots,i_s$ of~$A$
form an invertible matrix of size $s\times s$.
\end{enumerate}
\end{corollary}

Our next goal is to construct a bijection between the Gr\"obner fan
of~$I_L$ and the non-zero maximal minors of~$A$. The following
terminology will prove useful.

\begin{definition}\label{def-genGFan}
Let $J$ be an ideal in~$P$. 
\begin{enumerate}
\item[(a)] For a marked reduced Gr\"obner basis
$G = \{ (\LT_\sigma(g_1), g_1), \dots, (\LT_\sigma(g_k), g_k)\}$ 
of~$J$, we call $\LT_\sigma(G) = \{ \LT_\sigma(g_k),\dots, \LT_\sigma(g_k)\}$
the {\bf leading term set} of~$G$.

\item[(b)] The set $\LTGFan(J)$ of all leading term sets of marked
reduced Gr\"obner bases in $\GFan(J)$ is called the {\bf leading term Gr\"obner fan}
of~$J$. 
\end{enumerate}
\end{definition}

The following lemma provides some information about changing the basis of~$I_L$.
As above, by $A_{(i_1,\dots,i_s)}$ we denote the matrix
consisting of columns $i_1,\dots,i_s$ of a matrix~$A$.

\begin{lemma}\label{lem-easyassociated}
Let $L=(\ell_1, \dots, \ell_s)$ be a tuple of  $K$-linearly independent linear forms in~$P$,
let $I_L = \langle  L\rangle$, and let $A=(a_{ij})$ be the  coefficient matrix of~$L$. 
Moreover, let $L' = (\ell'_1, \dots, \ell'_r)$  be a further tuple 
of linear forms in~$I_L$, and let $A' \in\Mat_{r,n}(K)$ be its coefficient matrix.
\begin{enumerate}
\item[(a)]  The tuple~$L'$ is a minimal system of generators of~$I_L$ if and only if
$r=s$ and there exists a matrix $U \in \GL_s(K)$ such that $A' = U\cdot A$.

\item[(b)] A set of pairs $\{ (x_{i_1}, \ell'_1) \dots, (x_{i_r}, \ell'_r) \}$,
where $1\le i_1 < \cdots < i_r \le n$, 
is a marked reduced Gr\"obner basis of~$I_L$, if and only if $r=s$, the matrix $A_{(i_1,\dots,i_s)}$
is invertible, and $A' = (A_{(i_1, \dots, i_s)})^{-1}\cdot A$.
\end{enumerate}
\end{lemma}

\begin{proof} 
Claim~(a) follows from the fact that every tuple of 
minimal generators of~$I_L$ is also a basis of the $K$-vector space~$(I_L)_1$.

To prove~(b) we observe that a minimal Gr\"obner basis of a linear ideal 
is also a minimal set of generators of~$I_L$. This yields $r=s$. Moreover, it is reduced 
if and only if the submatrix $A'_{(i_1, \dots, i_s)}$ of~$A'$
is the identity matrix, and hence the conclusion follows from~(a).
\end{proof}

This lemma can also be interpreted in terms of the Pl\"ucker embedding of the
Gra{\ss}mannian ${\rm Gr}(s,n)$.
Now we are ready to present the key result for computing $\GFan(I_L)$.

\begin{theorem}{\bf (The Gr\"obner Fan of a Linear Ideal)}\label{thm-bijection}\\
Let $L=(\ell_1, \dots, \ell_s)$ be a tuple of $K$-linearly independent linear forms in~$P$,
let $I_L = \langle L\rangle$, and let $A=(a_{ij})$ be the  coefficient matrix of~$L$.  
Furthermore, let $M$ be the set of tuples $(i_1,\dots,i_s)$
such that $1\le i_1 < \cdots < i_s \le n$ and such that the corresponding
maximal minor of the matrix $A=(a_{ij})$ is non-zero.
\begin{enumerate}
\item[(a)] The map $\phi:\; \LTGFan(I_L) \longrightarrow M$ given by
$\phi(Z) = (i_1,\dots,i_s)$ for a tuple $Z = ( x_{i_1},\dots,x_{i_s} ) \in \LTGFan(I_L)$
with $1\le i_1 < \cdots < i_s\le n$ is well-defined and bijective.

\item[(b)] The  map $\psi:\; \GFan(I_L) \longrightarrow M$ given by
$\psi(G) = \phi(\LT_\sigma(G))$ for every $G\in \GFan(I_L)$ is well-defined and bijective.

\end{enumerate}
\end{theorem}

\begin{proof}
First we prove~(a).
To begin with, let us check that~$\phi$ is well-defined. For an element
$G\in \GFan(I_L)$, the tuple $Z=\LT_\sigma(G)$ satisfies $\langle Z\rangle
=\LT_\sigma(I_L)$. Let $Z = (x_{i_1},\dots,x_{i_s})$ with $1\le i_1< \cdots < i_s \le n$. 
Then Corollary~\ref{cor-CharLinZ}.c shows $\det(A_{(i_1,\dots,i_s)}) \ne 0$,
and therefore~$\phi$ is well-defined.

Since the map~$\phi$ is clearly injective, we still need to show that it is surjective.
Given $(i_1,\dots,i_s)\in M$, part~(b) of the lemma implies that
$(L')^{\rm tr} = A_{(i_1,\dots,i_s)}^{-1} \cdot L^{\rm tr}$ is a reduced Gr\"obner basis of~$I_L$
with leading term tuple $Z = (x_{i_1}, \dots, x_{i_s})$. 
Therefore we have $\phi(Z) = (i_1,\dots,i_s)$, and this proves the desired surjectivity.

To show~(b), we note that the map~$\psi$ is clearly well-defined and injective.
By the definition of $\LTGFan(I_L)$, it is also surjective. 
\end{proof}

Based on this theorem, we can compute the Gr\"obner fan of a linear ideal as follows.

\begin{corollary}{\bf (Computing the Gr\"obner Fan of a Linear Ideal)}\label{cor-computeGFlin}\\
Let $I_L=\langle \ell_1,\dots,\ell_s \rangle$ be an ideal in~$P$ generated
by linearly independent linear forms $\ell_1,\dots,\ell_s \in P_1$.
Then we can compute $\GFan(I_L)$ as follows.
\begin{enumerate}
\item[(1)] Let~$A$ be the coefficient matrix of~$L$, and let $S=\emptyset$.

\item[(2)] For every tuple $(i_1,\dots,i_s) \in M$, compute the maximal
minor $\det(A_{(i_1,\dots,i_s)})$ of~$A$.

\item[(3)] If $\det(A_{(i_1,\dots,i_s)}) \ne 0$, compute the 
vector $(L')^{\rm tr} = (A_{(i_1,\dots,i_s)})^{-1} \cdot L^{\rm tr}$ 
whose tuple of leading terms is $Z = ( x_{i_1}, \dots, x_{i_s} )$. 
Append the corresponding marked reduced Gr\"obner basis of~$I_L$ to~$S$.
Continue with the next tuple in Step~(2).

\item[(4)] Return the set $S=\GFan(I_L)$.
\end{enumerate}
\end{corollary}

Of course, depending on the numbers~$n$ and~$s$, computing $\binom{n}{s}$
determinants could be quite costly. The following remark may help us out.

\begin{remark}\label{rem-compGFan}
Let $I_L=\langle \ell_1,\dots,\ell_s \rangle$ be an ideal in~$P$ generated
by linearly independent linear forms $\ell_1,\dots,\ell_s \in P_1$ as above.

An alternative way of viewing the task to compute $\GFan(I_L)$ is
obtained by applying Corollary~\ref{cor-CharLinZ}. The complements
of the leading term sets of reduced Gr\"obner bases of~$I_L$ correspond uniquely
to sets of indeterminates whose residue classes form a $K$-basis of $P_1 / (I_L)_1$.
All sets of indeterminates whose residue classes are linearly independent 
in~$P_1/(I_L)_1$ are the independent sets of a linear matroid, and maximal
such sets are the bases of the matroid. The task of computing the bases of a linear
matroid has been studied intensively, and many algorithms are known,
see for instance the reverse search technique of D.\ Avis and K. Fukuda (cf.~\cite{AF}).
\end{remark}

Let us calculate the Gr\"obner fan of an explicit linear ideal.

\begin{example}\label{ex-restrictedGFan}
Let $P = \QQ[x, y, z, w]$, let $\ell_1=x+y -z +4w$, $\ell_2 =x-y-z$,
and let $I_L = \langle \ell_1, \ell_2\rangle$.
The set $\{\ell_1, \ell_2\}$ is a set of minimal generators of~$I_L$ 
and its coefficient matrix is
$$
A \;=\; \begin{pmatrix} 1 & \ \ 1 & -1 & 4 \cr  1 & -1 & -1& 0     
\end{pmatrix}.
$$
One $2\times2$-submatrix is singular. The others are 
$$
A_{12} = \left(\begin{smallmatrix} 1 & \;\; 1  \cr  1 & -1  \end{smallmatrix}\right),\,
A_{14} = \left(\begin{smallmatrix} 1 & 4  \cr  1 & 0  \end{smallmatrix}\right),\,
A_{23} = \left(\begin{smallmatrix} \;\; 1 & -1 \cr -1 & -1  \end{smallmatrix}\right),\,
A_{24} = \left(\begin{smallmatrix} \;\; 1 &  4 \cr -1 & 0  \end{smallmatrix}\right),\,
A_{34} = \left(\begin{smallmatrix} -1 & 4 \cr -1 & 0 \end{smallmatrix}\right) .
$$
Multiplying their inverses by~$A$ we get the matrices
$$
\Bigl( \begin{smallmatrix} 1 & 0 & -1 & 2 \cr 0 & \mathstrut 1 &\;\; 0 & 2 
   \end{smallmatrix}\Bigr) ,\,
\Bigl( \begin{smallmatrix} 1 & -1 & -1 & 0 \cr 0 & \;\; 1/2 & \;\; 0 & 1 
   \end{smallmatrix}\Bigr) ,\,
\Bigl( \begin{smallmatrix} \;\; 0 & 1 & 0 & \;\; 2 \cr -1 & \mathstrut 0 & 1 & -2  
   \end{smallmatrix}\Bigr) ,\,
\Bigl( \begin{smallmatrix} -1 & \; 1 & \;\; 1 & 0 \cr \;\; 1/2 & \; 0 & -1/2 & 1
   \end{smallmatrix}\Bigr) ,\,
\Bigl( \begin{smallmatrix} -1 & 1 & \; 1 & 0 \cr \;\; 0 & 1/2 & \; 0 & 1 
   \end{smallmatrix}\Bigr) . 
$$
They correspond to  the following marked reduced Gr\"obner bases of~$I_L$ 
which form the Gr\"obner fan of~$I_L$:
\begin{align*}
\{(x, x-z +2w), (y, y +2w) \}, & \quad \{(x, x -y - z), (w, w + \tfrac{1}{2}y)\}, \cr
\{(y, y+2w), (z, z-x-2w)\}, & \quad \{ (y, y-x+z), (w, w +\tfrac{1}{2}x -\tfrac{1}{2}z\}, \cr
\{ (z, z-x+y), (w,w +\tfrac{1}{2}y)\} .  \;\; &
\end{align*}
\end{example}

\bigskip\bigbreak
%%%%%%%%%%%%%%%%%%%%%%%%%%%%%%%%%%%%%%%%%%%%%%%%%%%%%%%%%%%%
%
% Section 4: Finding Z-Separating Re-embeddings
%
%%%%%%%%%%%%%%%%%%%%%%%%%%%%%%%%%%%%%%%%%%%%%%%%%%%%%%%%%%%%

\section{Finding Z-Separating Re-embeddings}
\label{Finding Z-Separating Re-embeddings}

In this section we show how to apply the Gr\"obner fan of the
linear part of an ideal to find tuples~$Z$ which are good candidates
for providing $Z$-separating re-embeddings of the ideal.
In~\cite{KLR2} we gave some answers to this question which use the computation 
of the Gr\"obner fan of the ideal~$I$ itself. 
Unfortunately, the computation of $\GFan(I)$ may be infeasible for large examples. 
The Gr\"obner fan of the ideal generated by the linear part of~$I$ is in general much smaller 
and thus provides a better set of candidate tuples~$Z$.

Using the definitions and notation introduced in the preceding sections, 
the following algorithm uses the Gr\"obner fan of the linear part of~$I$
in order to find $Z$-separating re-embeddings.

\begin{algorithm}{\bf ($Z$-Separating Re-embeddings 
via ${\boldsymbol \GFan(\langle \Lin_\M(I) \rangle)}$)}\label{alg-Emb}\\
Let $I\subseteq \M$ be an ideal of~$P$, and let $s \le \dim_K( \Lin_\M(I) )$.
Consider the following sequence of instructions.
\begin{enumerate}
\item[(1)] Using Corollary~\ref{cor-computeGFlin}, compute 
$\GFan(\langle \Lin_\M(I) \rangle)$.

\item[(2)] Form the set~$S$ of all tuples $Z=(z_1,\dots,z_s)$ such that
there is marked reduced Gr\"obner basis~$\overline{G}$ in 
$\GFan(\langle \Lin_\M(I) \rangle)$ for which $z_1,\dots,z_s$ are among the marked terms. 

\item[(3)] If $S = \emptyset$, return {\tt "No re-embedding found"}. 
While $S\ne \emptyset$, perform the following steps.

\item[(4)] Choose a tuple $Z = (z_1,\dots,z_s) \in S$ and remove it from~$S$.

\item[(5)] Using Remark~\ref{rem-checkZ}, check whether the ideal~$I$ is 
$Z$-separating. If it is, return~$Z$ and stop. Otherwise, 
continue with Step~(3).
\end{enumerate}
This is an algorithm which, if successful, finds a tuple of distinct indeterminates
$Z=(z_1,\dots,z_s)$ in~$X$ such that~$I$ is $Z$-separating. 

Moreover, if $s = \dim_K( \Lin_\M(I) )$ and the algorithm is successful then
the output tuple~$Z$ defines an optimal re-embedding of~$I$.
\end{algorithm}

\begin{proof}
Every tuple~$Z$ such that there exists a $Z$-separating re-em\-bed\-ding
of~$I$ is contained in the tuple of leading terms of a marked reduced Gr\"obner basis
of~$\Lin_{\M}(I)$ by Proposition~\ref{prop-cohSepPoly}.c.
The set of all possible such tuples~$Z$ is computed in Steps~(1) and~(2).
If the loop in Steps (3)-(5) finds a tuple~$Z$ such that~$I$ is $Z$-separating,
we are done.

In addition, if $s = \dim_K( \Lin_\M(I) )$ and the algorithm is successful,
then~\cite[Corollary~4.2]{KLR2}, shows that the $Z$-separating re-embedding 
of~$I$ is optimal.
\end{proof}

If we are looking for optimal re-embeddings and
use the method of Remark~\ref{rem-checkZ}.a to perform Step~(5), then
Algorithm~\ref{alg-Emb} is able to certify that no optimal \hbox{$Z$-se}\-parating 
re-embedding of~$I$ exists. However, we may have to compute some huge Gr\"obner bases. 
Moreover, the next remark points out some further limitations.

\begin{remark}\label{rem-onlysufficient}
Notice that Algorithm~\ref{alg-Emb} provides only a 
sufficient condition for detecting optimal re-embeddings of~$I$. 
On one side, it can happen that an optimal re-embedding is obtained using a subset 
of generators of a leading term ideal of $\Lin_\M(I)$ (see~\cite[Example~6.6]{KLR3}).
On the other side, it can happen that an optimal re-embedding cannot 
be achieved by a separating re-embedding, as shown in Example~\ref{ex-isotoK[x]}.
\end{remark}

The following example illustrates Algorithm~\ref{alg-Emb} at work.

\begin{example}\label{ex-searchingforoptembed}
Let $P=\QQ[x,y,z,w]$, let $F =(f_1, f_2, f_3)$, where $f_1 = x -y -w^2$, 
$f_2=x + y -z^2$, and $f_3 = z + w +z^3$, and let $I = \langle f_1,f_2,f_3\rangle$. 
\begin{enumerate}
\item[(1)] We obtain $\Lin_\M(I) = \langle z +w,\, x, \, y\rangle_K$ and the 
methods explained below return the two marked reduced 
Gr\"obner bases $\{ (x,x),\, (y,y),\, (z,z+w)\}$ and $\{ (x,x),\, (y,y),\, (w,w+z)\}$.

\item[(2)] We get $S =\{ (x, y, z),\, (x, y, w) \}$.

\item[(4)] We pick $Z = (x, y, z)$  and delete it from~$S$. 

\item[(5)] We construct an elimination ordering~$\sigma$ for~$Z$ and find that 
the minimal set of generators of~$\Lin_{\M}(\LT_\sigma(I))$ is $L = \{x\}$.
Therefore  $L\ne \{x,y,z\}$ and continue with the next iteration.

\item[(4)] Next we let $Z=(x,y,w)$ and let $S=\emptyset$.

\item[(5)] We construct an elimination ordering~$\sigma$ for~$Z$ 
and compute the minimal set of generators~$L$ of $\Lin_{\M}(\LT_\sigma(I))$. 
Since $L=\{x,y,w\}$, we return $Z = (x, y, w)$ and stop.
\end{enumerate}

To get  $Z$-separating polynomials, it suffices to 
replace $f_2$ with $f_2' = f_2-f_1$. 

To find the actual polynomials defining the optimal re-embedding, 
we compute the reduced $\sigma$-Gr\"obner basis of~$I$. It is
$$
( x -\tfrac{1}{2} z^6 -z^4 -z^2, \;  y + \tfrac{1}{2} z^6 +z^4,\;  w +z^3 +z).
$$
This tuple gives rise to a $\QQ$-algebra isomorphism $P/I\cong \QQ[z]$
via $x \mapsto \tfrac{1}{2} z^6 +z^4 +z^2$, $y \mapsto -\tfrac{1}{2} z^6 -z^4$, and
$w \mapsto -z^3 -z$.
\end{example}

\bigskip\bigbreak
%%%%%%%%%%%%%%%%%%%%%%%%%%%%%%%%%%%%%%%%%%%%%%%%%%%%%%%%%%%%%%%%%%%
%
% Section 5: Cotangent Equivalence Classes
%
%%%%%%%%%%%%%%%%%%%%%%%%%%%%%%%%%%%%%%%%%%%%%%%%%%%%%%%%%%%%%%%%%%%

\section{Cotangent Equivalence Classes}
\label{Cotangent Equivalence Classes}

The task to compute the Gr\"obner fan of $\langle\Lin_\M(I)\rangle$ in Algorithm~\ref{alg-Emb}
can be simplified when $\langle\Lin_\M(I)\rangle$ is a binomial linear ideal.
Recall that an ideal~$J$ in $P=K[x_1,\dots,x_n]$ is called a {\bf binomial ideal}
if it is generated by polynomials containing at most two terms in their support.

As before, we let $I = \langle f_1,\dots,f_r\rangle$
be an ideal in~$P$ with $f_i\in \M$. Letting $\ell_i= \Lin_\M(f_i)$ for $i=1,\dots,r$
and $L=(\ell_1,\dots,\ell_r)$, the linear part of~$I$ is $\Lin_\M(I) = \langle L \rangle_K$
and it generates the ideal $I_L = \langle \Lin_\M(I) \rangle$.

In the following we assume that the linear forms~$\ell_i$ are
{\bf binomials}, i.e., for $i=1,\dots, r$, we have 
$\ell_i = a_i\, x_{i_1} + b_i\, x_{i_2}$ with $a_i,b_i \in K$
and $i_1,i_2 \in \{1,\dots,n\}$. In this case the ideal $I_L$ is called a {\bf binomial 
linear ideal}.

Recall that, by Remark~\ref{rem-YgeneratesCot} and Corollary~\ref{cor-cohSepPoly}.c,
we can detect whether an ideal~$I$ is $Z$-separating by looking at the residue classes
of the entries of $Y = X \setminus Z$ in the cotangent space $\Cot_{\m}(R)=\m/\m^2$,
where $\m$ is the image of $\M = \langle x_1,\dots,x_n\rangle$ in $R=P/I$.
Moreover, we have $\Cot_\m(R) \cong P_1 / \langle L\rangle_K$.
This point of view leads us to the following definition.

\begin{definition}\label{def-cotequiv}
For every indeterminate $x_i\in X$, let $\bar{x}_i$ denote
its residue class in the cotangent space $\m / \m^2$ 
of $R=P/I$ at the origin.
\begin{enumerate}
\item[(a)] The relation~$\sim$ on~$X$ defined
by $x_i \sim x_j \Leftrightarrow \langle\bar{x}_i\rangle_K = \langle\bar{x}_j\rangle_K$
is an equivalence relation called {\bf cotangent equivalence}.

\item[(b)] An indeterminate $x_i\in X$ is called {\bf trivial}
if $\bar{x}_i=0$. The trivial indeterminates form the {\bf trivial
cotangent equivalence class} in~$X$.

\item[(c)] A non-trivial indeterminate $x_i \in X$ is called {\bf basic} if
its cotangent equivalence class consists only of~$x_i$. In this case, 
the cotangent equivalence class $\{ x_i \}$ is also called {\bf basic}.

\item[(d)] A non-trivial indeterminate $x_i \in X$ is called {\bf proper}
if its cotangent equivalence class contains at least two elements. In this case,
the cotangent equivalence class of~$x_i$ is also called {\bf proper}.
\end{enumerate}
\end{definition}

The meaning of these notions will become clear in the next theorems.
First we need a lemma which provides further information about the above definition.

\begin{lemma}\label{lem-equivclasses}
Let us assume to be in the above setting.
\begin{enumerate}
\item[(a)] The union~$U$ of the  supports of the elements in  
a minimal set of generators of the ideal~$I_L$ 
does not depend on the choice of a minimal set of generators. 

\item[(b)] The set of basic indeterminates is $X\setminus U$.

\item[(c)] The union of the sets of trivial and proper indeterminates is~$U$.
\end{enumerate}
\end{lemma}

\begin{proof}
To prove (a) we note that any minimal set of generators of the ideal
$\langle \Lin_\M(I) \rangle$ is also a minimal set of generators 
of the $K$-vector space $\Lin_\M(I)$.
Let~$A$ and~$B$ be two such sets. Since every linear form~$\ell$ in~$A$ 
is a linear combination of linear forms in~$B$, each indeterminate in $\Supp(\ell)$ 
is in the support of some linear form in~$B$. By interchanging the roles of~$A$ and~$B$, 
the conclusion follows.

Since $\sim$ is an equivalence relation on~$X$, to prove claims~(b) and~(c) it suffices 
to show that basic indeterminates are not in~$U$, while trivial and proper indeterminates
are in~$U$.
Firstly, let $x_i$ be a basic indeterminate. For a contradiction, assume that $x_i \in U$.
From $x_i\in U$ and the fact that $x_i$ is the only element in its equivalence class, 
we deduce that there is a polynomial in~$I$ of the form $x_i +q$ 
with $q \in \M^2$. Hence we get $\bar{x}_i = 0$, a contradiction to the fact that~$x_i$ is basic.
Secondly, let $x_i$ be trivial.  Then there is a polynomial in~$I$ of the form $x_i+q$
with $q \in \M^2$, and hence we get $x_i \in U$. 
Thirdly, let $x_i$ be proper. Then there exist another indeterminate $x_j$ and a polynomial in~$I$ 
of the form  $a x_i + b x_j+q$ with $a,b\in K\setminus \{0\}$ and $q \in \M^2$. Hence we get $x_i \in U$.
\end{proof}

Our next step is to order the indeterminates in the cotangent equivalence classes
using a term ordering. The following definition will come in handy.

\begin{definition}\label{def-leadsetofPr}
Let $E= \{x_{i_1}, \dots, x_{i_p} \}$ be a proper equivalence class in~$X$,
and let $\sigma$ be a term ordering on ${\mathbb T}^n$
with $x_{i_1} >_\sigma \cdots  >_\sigma x_{i_p} $. Then the set
$E\setminus\{ x_{i_p} \}$ is called the {\bf $\sigma$-leading set} of~$E$ and
denoted by~$E^\sigma$.
\end{definition}

\begin{notation}
Let $\sigma$ be a term ordering on~${\mathbb T}^n$.
In accordance with the notation introduced in Proposition~\ref{prop-cohSepPoly},
the unique minimal set of indeterminates generating the ideal 
$\LT_\sigma(I_L)$ will be denoted by~$S_\sigma$.
\end{notation}

In the following theorem we give an explicit representation of~$S_\sigma$
for every term ordering~$\sigma$ on $\mathbb{T}^n$
and show the importance of $\sigma$-leading sets.
This implies  a description of $\LTGFan(I_L)$
which has several advantages compared to the general 
description given in Theorem~\ref{thm-bijection}. 
For instance, it does not have to deal with huge matrices.

\begin{theorem}\label{thm-shapeofSsigma}
Let $E_0$ be the trivial equivalence class, and let $E_1, \dots, E_q$ be
the proper equivalence classes in~$X$. Let $I_L = \langle \Lin_\M(I)\rangle$
be the ideal generated by the linear parts of the polynomials in~$I$.
\begin{enumerate}
\item [(a)] Let~$\sigma$ be a term ordering on $\mathbb{T}^n$.
Then we have $S_\sigma = E_0 \cup E_1^\sigma \cup \cdots \cup E_q^\sigma$,
and hence $\# S_\sigma  = \# E_0 + \sum_{i=1}^q \# E_i - q$.

\item[(b)] For $i=1,\dots,q$, let $E_i^\ast$ be a set obtained from~$E_i$ 
by deleting one of its elements. Then there exists a term ordering~$\sigma$
on $\mathbb{T}^n$ such that we have
$S_\sigma = E_0 \cup E_1^\ast \cup \cdots \cup E_q^\ast$.

\item[(c)]  Let $\Sigma$ be the set of all sets of the form 
$E_0 \cup E_1^\ast \cup \cdots \cup E_q^\ast$, where~$E_i^\ast$ is obtained
from the set~$E_i$ by deleting one of its elements.
Then the map $\phi: \LTGFan(I_L)  \longrightarrow \Sigma$ given by 
$\phi(S_\sigma) = E_0 \cup E_1^\sigma \cup \cdots \cup E_q^\sigma$ 
is well-defined and bijective.

\item[(d)]  We have
$\#\LTGFan(I_L) = \#\GFan(I_L)= \prod\limits_{i=1}^q \# E_i$.
\end{enumerate}
\end{theorem}

\begin{proof} 
To prove claim~(a) we observe that the inclusion $E_0 \subseteq S_\sigma$ follows from 
$E_0 \subseteq \Lin_\M(I)$, and that  the inclusion 
$S_\sigma \subseteq E_0 \cup E_1\cup \cdots \cup E_q$
follows from Lemma~\ref{lem-equivclasses}.c.

For $k\in \{1,\dots,q\}$, we write the proper equivalence class 
$E_k= \{ x_{i_1}, \dots, x_{i_p} \}$ 
such that $x_{i_1} >_\sigma \cdots  >_\sigma x_{i_p}$. 
Using the definition of a proper equivalence class, it follows that 
$x_{i_m} - x_{i_p} \in \Lin_\M(I)$, so that
$\LT_\sigma(x_{i_m} - x_{i_p}) = x_{i_m}$ 
for every $m \in \{1,\dots,p-1\}$.
Hence we have proved 
$$
E_0 \cup E_1^\sigma \cup \cdots \cup E_q^\sigma \;\subseteq\;  
S_\sigma \;\subseteq\;   E_0 \cup E_1\cup \cdots \cup E_q .
$$ 
As the $\sigma$-smallest element in each proper equivalence class does not belong to~$S_\sigma$, 
claim~(a) follows.
 
Claim~(b) follows from~(a) if we show that there exists a term ordering~$\sigma$ 
such that $E_i^\sigma = E_i^\ast$ for $i=1,\dots,q$.
By definition, the sets~$E_i$ are pairwise disjoint. Consequently, a term ordering 
which solves the problem can be chosen as a block term ordering, and hence it 
suffices to consider the case $q=1$. 
So, let $E_1 = \{x_{i_1}, \dots, x_{i_p} \}$. W.l.o.g.\ assume 
that $E_1^\ast = \{x_{i_2}, \dots, x_{i_p} \}$. As we observed before, we have 
$x_{i_k} - x_{i_1} \in \Lin_\M(I)$ for $k=2,\dots,p$.
To finish the proof it suffices to take a term ordering~$\sigma$ such that 
$x_{i_k} >_\sigma x_{i_1}$ for $k=2,\dots,p$.
  
Claim~(c) follows immediately from~(b).
To prove claim~(d) we note that the first equality is obvious.  
To show $\# \LTGFan(I_L) =  \prod_{i=1}^q \# E_i$, it suffices to 
deduce from~(b) that the number of the leading term ideals of~$I_L$ 
equals the number of $q$-tuples of indeterminates $x_i$, exactly one chosen in
each proper equivalence class.
\end{proof}

Now we are ready to classify the indeterminates in~$X$ which can be used for a 
$Z$-separating re-embedding of~$I$ as follows.

\begin{theorem}\label{thm-classifycij}
Let~$Z$ be a tuple of indeterminates from~$X$ such that there exists a
$Z$-separating re-embedding of~$I$, and let $Y = X \setminus Z$.
\begin{enumerate}
\item[(a)] The basic indeterminates of~$X$ are contained in~$Y$.

\item[(b)] Each proper equivalence class in~$X$ contains at least one element of~$Y$.

\item[(c)] If $\# Z = \dim_K(\Lin_\M(I))$, then the $Z$-separating 
re-em\-bed\-ding of~$I$ is optimal, the trivial indeterminates of~$X$ are contained 
in~$Z$, and each proper equivalence class in~$X$ contains exactly one element of~$Y$.

\end{enumerate}
\end{theorem}

\begin{proof}
To prove~(a), let~$\sigma$ be a $Z$-separating term ordering for~$I$,
let~$S_\sigma$ be the minimal set of indeterminates generating $\LT_\sigma(I_L)$,
and let~$U$ be the union of the supports of the elements in a minimal set of generators 
of~$I_L$. Note that Proposition~\ref{prop-cohSepPoly}.c 
implies $Y \supseteq X \setminus S_\sigma$. From $S_\sigma \subseteq X$ we deduce the inclusion 
$Y \supseteq X \setminus U$, and thus the claim follows from Lemma~\ref{lem-equivclasses}.b.
 
Claim~(b) follows from the mentioned relation $Y \supseteq X \setminus S_\sigma$ and 
Theorem~\ref{thm-shapeofSsigma}.a.

Finally, we prove (c). If $\#(Z) = \dim _K(\Lin_\M(I))$ then Corollary~\ref{cor-cohSepPoly}.a implies that 
the $Z$-separating re-embedding of~$I$ is optimal. Moreover,
Corollary~\ref{cor-cohSepPoly}.b implies that $Z = S_\sigma$, and hence the claim
follows from Theorem~\ref{thm-shapeofSsigma}.a.
\end{proof}

Based on the preceding results and on Algorithm~3.8 in~\cite{KSL} for
computing the cotangent equivalence classes, we can now check effectively
whether a given ideal~$I$ admits a $Z$-separating embedding. 
Notice that we are excluding some trivial cases (namely $n=1$ and $I=\M$) 
in order to be able to apply \cite[Algorithm~3.8]{KSL}, but these
cases can be dealt with easily by a direct computation.

\begin{algorithm}{\bf ($Z$-sep.\ Re-embeddings Via Cotangent Equivalence)}\label{alg-compEmb}\\
Let $I\subsetneq \M$ be an ideal in $P=K[x_1,\dots,x_n]$, where $n\ge 2$, and let
$X=(x_1,\dots,x_n)$. Consider the following sequence of instructions.
\begin{enumerate}
\item[(1)] Compute the trivial cotangent equivalence class~$E_0$
and also the proper cotangent equivalence classes $E_1,\dots,E_q$.

\item[(2)] Let $S=\emptyset$.

\item[(3)] Turn each set $Z_0\cup Z_1\cup \cdots \cup Z_q$ such that 
$Z_0 \subseteq E_0$ and $Z_i \subsetneq E_i$ for $i=1,\dots,q$ into a tuple~$Z$ 
and perform the following steps. 

\item[(4)] Using Remark~\ref{rem-checkZ}, check whether the
ideal~$I$ is $Z$-separating. If it is, append~$Z$ to~$S$.

\item[(5)]  Continue with Step~(3) using the next tuple $Z$ until 
all tuples have been dealt with. Then return~$S$ and stop.
\end{enumerate}

Then the following two claims hold.
\begin{enumerate}
\item[(a)] This is an algorithm which computes the set~$S$ of all tuples~$Z$ 
of distinct indeterminates in~$X$ such that there exists a $Z$-sep\-a\-ra\-ting 
re-embedding of~$I$.

\item[(b)] Assume that Step~(3) is replaced by the following step. 

\begin{enumerate}
\item [(3')] Turn each set $E_0\cup E_1^\ast\cup \cdots \cup E_q^\ast$, 
where $E_i^\ast$ is obtained from~$E_i$ by deleting one element, 
into a tuple~$Z$ and perform the following steps. 
\end{enumerate}
Then the result is an algorithm which computes the set~$S$ of all tuples~$Z$ 
of distinct indeterminates in~$X$ such that there exists an optimal
$Z$-separating re-embedding of~$I$.
\end{enumerate}
\end{algorithm}

\begin{proof}
Both claims follow from Theorem~\ref{thm-shapeofSsigma} and Theorem~\ref{thm-classifycij}.  
\end{proof}

For an example to illustrate this algorithm, we refer the reader to the next section.

\bigskip\bigbreak
%%%%%%%%%%%%%%%%%%%%%%%%%%%%%%%%%%%%%%%%%%%%%%%%%%%%%%%%%%%%%%%%%%
%
% Section 6: Application to Border Basis Schemes
%
%%%%%%%%%%%%%%%%%%%%%%%%%%%%%%%%%%%%%%%%%%%%%%%%%%%%%%%%%%%%%%%%%%

\section{Application to Border Basis Schemes}
\label{Application to Border Basis Schemes}

In this section we apply the methods developed above to the ideals defining border basis schemes. 
These affine schemes are moduli spaces of 0-dimensional ideals which are canonically embedded
into very high dimensional affine spaces. To study them carefully, it is imperative
to re-embed them into lower dimensional affine spaces. For instance, one important question
is whether a given border basis scheme is an {\bf affine cell}, i.e., isomorphic
to an affine space. As we shall recall below, the natural generators of the defining
ideals of border basis schemes have binomial linear parts, so that the theory developed
in the preceding section is perfectly suited to re-embed border basis schemes. In fact, 
we also provide some information complementing Theorem~\ref{thm-classifycij} in this situation.

In the following we assume that the reader has a basic knowledge of border basis
theory, e.g., to the extent it is covered in~\cite[Section~6.4]{KR2}.
Let $\OO = \{t_1,\dots,t_\mu\}$ be an order ideal of terms in~$\mathbb{T}^n$,
and let $\partial\OO = \{ b_1,\dots,b_\nu\}$ be its border. By replacing the coefficients
of an $\OO$-border prebasis with indeterminates $c_{ij}$, we obtain the {\bf generic
$\OO$-border prebasis} $G=\{g_1,\dots,g_\nu\}$, where
$g_j = b_j - \sum_{i=1}^\mu \, c_{ij}\, t_i$
for $i=1,\dots,\mu$ and $j=1,\dots,\nu$. Furthermore, let~$C$ be the set of indeterminates
$C= \{ c_{ij} \mid i=1,\dots,\mu;\; j=1,\dots,\nu\}$. The ideal $I(\BO)$ in $K[C]$
defining the border basis scheme $\BO$ can be constructed in several ways (see~\cite{KR3},
\cite{KR4}, \cite{KLR1}):
\begin{enumerate}
\item[(1)] As in~\cite[Definition~3.1]{KR3}, construct the {\bf generic multiplication 
matrices} $\mathcal{A}_1,\dots, \mathcal{A}_n \in \Mat_\mu(K[C])$ and let $I(\BO)$ be the ideal 
in~$K[C]$ generated by all entries of all commutators $\mathcal{A}_i\, \mathcal{A}_j - 
\mathcal{A}_j\, \mathcal{A}_i$ with $1\le i<j\le n$.

\item[(2)] Construct the set of next-door generators $\ND_\OO$
and the set of across-the-rim generators $\AR_\OO$ of~$I(\BO)$ and take
the union $\ND_\OO \cup \AR_\OO$ (see below).
\end{enumerate}

For us, the most important properties of these sets of polynomials $\{f_1,\dots,f_r\}$
are that~$f_i$ consists for $i=1,\dots,r$ of a linear and a quadratic part. In the following 
we describe these homogeneous components in more detail. We begin by making the
construction in~(2) explicit.

\begin{definition}{\bf (Neighbour Generators)}\label{def-neighbours}\\
Let $\OO = \{t_1,\dots,t_\mu \}$ be an order ideal in~$\mathbb{T}^n$ with border
$\partial\OO = \{ b_1,\dots,b_\nu \}$.
Let $\mathcal{A}_1,\dots,\mathcal{A}_n$ be the
generic multiplication matrices, and for $j=1,\dots,\mu$
let $c_j=(c_{1j},\dots,c_{\mu j})\tr$ be the $j$-th column
of~$(c_{ij})$.
\begin{enumerate}
\item[(a)] Let $j, j'\in\{1,\dots,\nu\}$ be such that
$b_j=x_\ell b_{j'}$ for some $\ell\in\{1,\dots,n\}$. Then $b_j,b_{j'}$
are called {\bf next-door neighbours} and the tuple of polynomials
$(c_j - \mathcal{A}_\ell c_{j'})\tr$ is denoted by $\ND(j, j')$.

\item[(b)] The union of all entries of the tuples $\ND(j, j')$ is called the set of
{\bf next-door generators} of~$I(\BO)$ and is denoted by
$\ND_\OO$.

\item[(c)] Let $j, j'\in\{1,\dots,\nu\}$ be such that
$b_j=x_\ell b_m$ and $b_{j'}=x_k b_m$ for some $m\in \{1,\dots,\nu\}$.
Then $b_j,b_{j'}$ are called {\bf across-the-corner neighbours}.

\item[(d)] Let $j, j'\in\{1,\dots,\nu\}$ be such that $b_j = x_\ell t_m$ and
$b_{j'} = x_k t_m$ for some $m\in\{1,\dots,\mu\}$. Then $b_j,b_{j'}$
are called {\bf across-the-rim neighbours} and the tuple of polynomials
$(\mathcal{A}_k c_j - \mathcal{A}_\ell c_{j'})\tr$
is denoted by $\AR(j, j')$.

\item[(e)] The union of all entries of the tuples $\AR(j, j')$ is called the set
of {\bf across-the-rim generators} of~$I(\BO)$ and is denoted by $\AR_\OO$.

\item[(f)] The polynomials in $\ND_\OO \cup \AR_\OO$ are called the 
{\bf neighbour generators} of~$I(\BO)$.
\end{enumerate}
\end{definition}

In~\cite[Proposition~4.1]{KR3}, it is shown that the polynomials corresponding
to across-the-corner neighbours are not necessary to generate $I(\BO)$
and that the neighbour generators are precisely the non-trivial 
entries of the commutators $\mathcal{A}_i\, \mathcal{A}_j - \mathcal{A}_j\, 
\mathcal{A}_i$.

An important property of the neighbour generators is that they are homogeneous with respect
to the following multigrading. Recall that the logarithm of a term $t=x_1^{\alpha_1} \cdots
x_n^{\alpha_n}$ is defined by $\log(t) = (\alpha_1,\dots,\alpha_n)$.

\begin{definition}
The $\mathbb{Z}^n$-grading on~$K[C]$ defined by $\deg_W(c_{ij}) = \log(b_j)-\log(t_i)$
for $i=1,\dots,\mu$ and $j=1,\dots,\nu$ is called the {\bf arrow grading}. 
\end{definition}

As mentioned above, the neighbour generators of~$I(\BO)$ have (standard) degree two and no 
constant term. Their linear parts can be described in detail as follows.

\begin{proposition}{\bf (Linear Parts of Neighbour Polynomials)}\label{prop-LinParts}\\
Let $\OO=\{t_1,\dots,t_\mu\}$ be an order ideal in~$\mathbb{T}^n$
with border $\partial\OO = \{b_1,\dots, b_\nu\}$.
\begin{enumerate}
\item[(a)] Let $j, j'\in\{1,\dots,\nu\}$ be such that $b_j,b_{j'}$
are next-door neighbours, i.e., such that $b_j = x_\ell b_{j'}$, 
and let $i\in \{1,\dots,\mu\}$.
Then the linear part of the corresponding next-door generator
$ND(j,j')_i$ is (up to sign) given by
$$
\begin{cases}
c_{ij}-c_{i' j'}  &  \text{ if $x_\ell$ divides $t_i$,} \\
c_{ij} & \text{ otherwise.}
\end{cases}
$$
The polynomial $f=\ND(j,j')_i$ is homogeneous with respect to the arrow degree
with $\deg_W(f)=\deg_W(c_{ij})$.

\item[(b)]  Let $j,j'\in \{1,\dots,\nu\}$ be such that $b_j, b_{j'}$
are across-the-rim neighbours, i.e., such that $b_j = x_\ell t_{m'}$ 
and $b_{j'} = x_k t_{m'}$ for some $m'\in\{1,\dots,\mu\}$. For $m\in\{1,\dots,\mu\}$, 
the non-zero linear part of the corresponding across-the-rim generator $\AR(j,j')_m$
is (up to sign) given by
$$
\begin{cases} c_{ij}-c_{i'j'} & \text{ if $t_m= x_k t_i= x_\ell t_{i'}$,}\\
c_{ij} & \text{ if $t_m = x_k t_i$, but $x_\ell$ does not divide~$t_m$,}\\
c_{i'j'} & \text{ if $t_m = x_\ell t_{i'}$, but $x_k$ does not divide~$t_m$.}
\end{cases}
$$ 
The polynomial $g=\AR(j,j')_m$ is homogeneous with respect to the arrow degree 
with $\deg_W(g)=\deg_W(c_{mj}) + e_k = \deg(c_{mj'}) + e_\ell$.

\end{enumerate}
\end{proposition}

\begin{proof} 
The claims for the linear parts are shown in~\cite[Corollary~2.8]{KSL}.
The claims for the arrow degree follow from~\cite[Lemma~3.4]{KSL}, and the description
of the linear parts. 
\end{proof}

To describe the quadratic parts of the neighbour generators in detail, 
the following concepts are convenient.

\begin{definition}
Let $\OO=\{t_1,\dots,t_\mu\}$ be an order ideal in~$\mathbb{T}^n$.
\begin{enumerate}
\item[(a)] The {\bf rim} $\OO^\nabla$ of~$\OO$ consists of all terms~$t_i$
such that $x_k t_i \in \partial\OO$ for some $k \in \{1,\dots,n\}$.
The indeterminate $c_{ij}$ is called a {\bf rim indeterminate}
if $t_i \in \OO^\nabla$, and the set of all rim 
indeterminates is denoted by $C^\nabla$.

\item[(b)] The {\bf interior} of~$\OO$ is $\OO^\circle = \OO \setminus \OO^\nabla$.
The indeterminate $c_{ij}$ is called an {\bf interior indeterminate}
if $t_i\in \OO^\circle$, and the set of all interior indeterminates is denoted by~$C^\circle$.
\end{enumerate}
\end{definition}

Clearly, we have a disjoint union $C = C^\nabla \cupdot C^\circle$.
One more definition, and we are ready to go. The following notion 
was introduced in~\cite[Section~4.1]{Hui1}.

\begin{definition}\label{def-exposedterm}
Let $j\in\{1,\dots,\nu\}$, and let $\ell\in\{1,\dots,n\}$.
Then the border term~$b_j$ is called {\bf $x_\ell$-exposed}
if it is of the form $b_j = x_\ell \, t_i$ with
$i \in \{1,\dots,\mu\}$. In this case we also say that~$t_i$
{\bf $x_\ell$-exposes} the border term~$b_j$.
\end{definition}

Finally we are ready to describe the homogeneous components of (standard)
degree two of the neighbour generators of $I(\BO)$.

\begin{proposition}{\bf (Quadratic Parts of Neighbour Generators)}\label{prop-QuadrParts}\\
Let $\OO=\{t_1,\dots,t_\mu\}$ be an order ideal in~$\mathbb{T}^n$
with border $\partial\OO = \{b_1,\dots, b_\nu\}$.
\begin{enumerate}
\item[(a)] Let $j, j'\in\{1,\dots,\nu\}$ be such that $b_j,b_{j'}$
are next-door neighbours with $b_j = x_\ell b_{j'}$, and let $i\in \{1,\dots,\mu\}$.
Assume that $b_{\lambda_1}, \dots, b_{\lambda_s} \in\partial\OO$ are the
$x_\ell$-exposed border terms, and let $b_{\lambda_p} = x_\ell \, t_{\rho_p}$
for $p=1,\dots,s$. Then the quadratic terms in the support of~$\ND(j,j')_i$
are the products $c_{i \lambda_p}\, c_{\rho_p j'}$ with $p=1,\dots,s$.

In particular, all terms in the quadratic part are of the form $c_{i\lambda} c_{\rho j'}$
with $\lambda\in \{1,\dots,\nu\}$ and a rim indeterminate $c_{\rho j'}$.

\item[(b)] Let $j,j'\in \{1,\dots,\nu\}$ be such that $b_j, b_{j'}$
are across-the-rim neighbours with $b_j = x_\ell t_{m'}$ and $b_{j'} = x_k t_{m'}$
for some $m'\in\{1,\dots,\mu\}$. For every $m\in\{1,\dots,\mu\}$,
the quadratic terms in the support of $\AR(j,j')_m$ are of the following types.
\begin{itemize}
\item[(1)] Let $b_{\kappa_1},\dots, b_{\kappa_s} \in \partial\OO$ be the
$x_k$-exposed border terms, and let us write $b_{\kappa_p} = x_\ell\, t_{\rho_p}$
for $p=1,\dots,s$. 
Then the terms  $c_{m \kappa_p} \, c_{\rho_p j}$ may appear in the support of $\AR(j,j')_m$.

\item[(2)] Let $b_{\lambda_1}, \dots,b_{\lambda_u} \in \partial\OO$ be the
$x_\ell$-exposed border terms, and let us write $b_{\lambda_q} = x_\ell t_{\sigma_q}$
for $q=1,\dots,u$. Then the terms  $c_{m \lambda_q}\, c_{\sigma_q j'}$
may appear in the support of $\AR(j,j')_m$.
\end{itemize}

In particular, all terms in the quadratic part
are of the form $c_{m\lambda} c_{\rho j}$ or $c_{m\lambda} c_{\rho j'}$
with $\lambda \in \{1,\dots,\nu\}$ and rim indeterminates $c_{\rho j}$ and 
$c_{\rho j'}$, respectively. 

\end{enumerate}
\end{proposition}

\begin{proof}
Let $G=\{g_1,\dots,g_\nu\}$ be the generic $\OO$-border prebasis. We shall use the construction
of the neighbour generators using the lifting of neighbour syzygies (see~\cite[Section~5]{KK}
and \cite[Proposition~6.4.34]{KR2}).

First we prove~(a). The polynomials in $\ND(j,j')$ are the coefficients of 
$t_1,\dots,t_\mu$ in the reduction of $x_\ell g_{j'} - g_j$, viewed as a polynomial in
$K[C][x_1,\dots,x_n]$. The leading terms $x_\ell b_{j'}$ and $b_j$ cancel by
definition. Hence we only have to reduce the $x_\ell$-exposed border terms
in $x_\ell g_{j'}$. The coefficient of~$b_{\lambda_p}$ in $x_\ell g_{j'}$
is $c_{\rho_p j'}$. The coefficient of~$t_i$ in~$g_{\lambda_p}$ is $c_{i \lambda_p}$.
Therefore the coefficient of~$t_i$ in $c_{\rho_p j'}\, g_{\lambda_p}$
is $c_{\rho_p j'} c_{i \lambda_p}$ and it will be a part of the coefficient
of~$t_i$ in the final result of the reduction, i.e., in $\ND(j,j')_i$.

Next we prove~(b). Again the entries of $\AR(j,j')$ are obtained
as the coefficients of $t_1,\dots,t_\mu$ in the reduction of $x_k g_j - x_\ell g_{j'}$,
and again the highest terms $x_k b_j$ and $x_\ell b_{j'}$ cancel.
To reduce $x_k g_j$, we have to reduce the $x_k$-exposed border terms
$b_{\kappa_1},\dots, b_{\kappa_s}$. The coefficient of $b_{\kappa_p}$ in $x_k g_j$ is
$c_{\rho_p j}$. The polynomial $\AR(j,j')_m$ has arrow degree $\deg_W(c_{mj})+e_k$.
It is therefore the coefficient of~$t_m$ in the reduction of $x_k g_j - x_\ell
g_{j'}$. The coefficient of $t_m$ in $g_{\kappa_p}$ is~$c_{m \kappa_p}$.
Hence the coefficient of~$t_m$ in $c_{\rho_p j}\, g_{\kappa_p}$ is
$c_{\rho_p j} c_{m \kappa_p}$, and it may appear in $\AR(j,j')_m$.
These quadratic terms are the ones listed in case~(1).

The analysis of the terms in the reduction of $x_\ell g_{j'}$
is completely analogous and leads to the quadratic terms in case~(2).
\end{proof}

Based on this detailed study of the support of the neighbour generators,
we are now able to provide some additional information on the distribution of
the rim indeterminates in the cotangent equivalence classes.

\begin{theorem}\label{thm-RimIndets}
Let~$Z$ be a tuple of indeterminates from~$C$ such that there exists a
$Z$-separating re-embedding of~$I(\BO)$ and let $Y = C \setminus Z$.
\begin{enumerate}
\item[(a)] All basic indeterminates of~$C$ are rim indeterminates.

\item[(b)] Each proper cotangent equivalence class in~$C$ contains a rim indeterminate.
\end{enumerate}
\end{theorem}

\begin{proof}
For the proof of~(a), we assume that $c_{ij} \in C$ is a basic indeterminate and that
$t_i\in \OO^\circle$. If there exists an indeterminate $x_\ell$ such that
$x_\ell b_j \in \partial\OO$ then we let $j'\in\{1,\dots,\nu\}$ such that 
$b_{j'} = x_\ell b_j$. As we have $t_i\in\OO^\circle$, there exists an index
$i'\in\{1,\dots,\mu\}$ such that $t_{i'} = x_\ell t_i \in \OO$.
Hence we have a next-door neighbour pair $(b_j,b_{j'})$ and 
Proposition~\ref{prop-LinParts} implies $c_{ij} \sim c_{i' j'}$. This contradicts
the hypothesis that $c_{ij}$ is a basic indeterminate. Thus no next-door
neighbour pair $(b_j,b_{j'})$ exists.

Since the border~$\partial\OO$ is connected with respect to the neighbour
relations (see~\cite[Proposition~19]{KK}), there exists 
an across-the-rim neighbour pair $(b_j,b_{j'})$.
Thus we may assume that there are $m'\in\{1,\dots,\mu\}$ and $k,\ell \in\{1,\dots,n\}$
such that $b_j = x_\ell t_{m'}$ and $b_{j'} = x_k t_{m'} \in\partial\OO$
for some $j'\in \{1,\dots,\nu\}$. Now, if $t_i$ is not divisible by~$x_\ell$,
then Proposition~\ref{prop-QuadrParts} implies that $c_{ij}$ is a trivial indeterminate,
in contradiction to the hypothesis. Hence $t_i$ has to be divisible by~$x_\ell$.
The term $t_i x_k/x_\ell$ cannot be in the border of~$\OO$, because then also
$t_i x_k$ would be outside~$\OO$, in contradiction to the hypothesis $t_i\in\OO^\circle$.
Therefore there exists an index $i'\in\{1,\dots,\mu\}$ such that $x_\ell t_{i'} = x_k t_i$.
By Proposition~\ref{prop-LinParts}, we get $c_{ij} \sim c_{i'j'}$, in contradiction
to the hypothesis that $c_{ij}$ is basic.

Now we show~(b). For a contradiction, assume that $c_{ij}\in C$ is proper
and that every element in the cotangent equivalence class of~$c_{ij}$
is an interior indeterminate. Since~$\OO$ is an order ideal, the term~$t_i$
is not a multiple of $b_j\notin \OO$. Thus there exists
an index $\ell\in\{1,\dots,n\}$ such that the arrow degree $\deg_W(c_{ij}) = \log(b_j) - \log(t_i)$
of~$c_{ij}$ has a positive $\ell$-th component. Let $j'\in\{1,\dots,\nu\}$ be such that
the $\ell$-th component of $\log(b_{j'})$ is zero. (For instance, let $k\ne \ell$
and consider the unique term of the form $b_{j'} = x_k^N \in\partial\OO$ with 
$N\ge 1$.) Since the border is connected, we can find a sequence of border terms
$b_j = b_{j_1} \sim \cdots \sim b_{j_q} = b_{j'}$ such that $b_{j_p}$ and $b_{j_{p+1}}$
are next-door or across-the-rim neighbours for $p=1,\dots,q-1$.

Notice that it is not possible to find indeterminates $c_{i_1} \sim \cdots\sim c_{i_q j_q}$
because the equality of arrow degrees $\deg_W(c_{i_q j_q}) = \deg_W(c_{i_1 j_1})$ 
would imply that the $x_\ell$-degree of $t_{i_q}$ is negative. Let us try to construct
such a sequence of cotangent equivalences inductively. When we have found terms 
$t_{i_1},\dots,t_{i_p}\in\OO$ with the property that $c_{i_1 j_1} \sim \cdots\sim c_{i_p j_p}$
and try to find a term $t_{i_{p+1}}\in \OO$ such that 
$c_{i_p j_p} \sim c_{i_{p+1} j_{p+1}}$, three things can happen:
\begin{enumerate}
\item[(1)] A term $t_{i_{p+1}}$ of the desired kind exists in~$\OO^\circle$.
\item[(2)] A term $t_{i_{p+1}}$ of the desired kind exists in~$\OO^\nabla$.
\item[(3)] No term $t_{i_{p+1}}$ of the desired kind exists in~$\OO$ because
one of the components of $\log(t_{i_{p+1}})$ would be negative.
\end{enumerate}
In case~(1), we can continue our inductive construction for one further step.
By the hypothesis that $c_{ij}$ is not equivalent to a rim indeterminate,
case~(2) never occurs. Hence case~(3) has to happen for some $p\in\{ 1,\dots,q-1\}$.
In this case we have $\bar{c}_{i_p j_p} = 0$ by Proposition~\ref{prop-LinParts},
and thus $\bar{c}_{ij} = 0$. Hence we have arrived at a contradiction to the hypothesis 
that~$c_{ij}$ is proper, and the proof is complete.
\end{proof}

The following example illustrates the results of this section and the preceding one.
The readers may also check that it verifies some claims in~\cite[Remark~7.5.3]{Hui1}.

\begin{example}\label{ex-Section6}
In $P =\QQ[x,y]$, consider the order ideal
$\OO = \{t_1,\dots,t_8\}$ given by $t_1=1$, $t_2=y$, $t_3=x$, $t_4=y^2$, $t_5=xy$, $t_6=x^2$, 
$t_7=y^3$, and $t_8 =xy^2$. Then we have $\partial\OO = \{b_1,\dots,b_5\}$ with $b_1=x^2y$, 
$b_2=x^3$, $b_3=y^4$, $b_4 = xy^3$, and $b_5=x^2y^2$. 
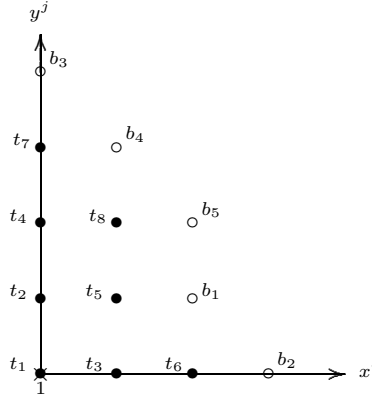
\begin{center}
	\makebox{\beginpicture
		\setcoordinatesystem units <1cm,1cm>
		\setplotarea x from 0 to 4, y from 0 to 4.5
		\axis left /
		\axis bottom /
		\arrow <2mm> [.2,.67] from  3.5 0  to 4 0
		\arrow <2mm> [.2,.67] from  0 4  to 0 4.5
		\put {$\scriptstyle x^i$} [lt] <0.5mm,0.8mm> at 4.1 0.1
		\put {$\scriptstyle y^j$} [rb] <1.7mm,0.7mm> at 0 4.6
		\put {$\bullet$} at 0 0
		\put {$\bullet$} at 1 0
		\put {$\bullet$} at 0 1
		\put {$\bullet$} at 1 1
	         \put {$\bullet$} at 1 2
		\put {$\bullet$} at 0 2
		\put {$\bullet$} at 0 3
		\put {$\bullet$} at 2 0
		\put {$\scriptstyle 1$} [lt] <-1mm,-1mm> at 0 0
		\put {$\scriptstyle t_1$} [rb] <-1.3mm,0.4mm> at 0 0
		\put {$\scriptstyle t_3$} [rb] <-1.3mm,0.4mm> at 1 0
		\put {$\scriptstyle t_2$} [rb] <-1.3mm,0mm> at 0 1
		\put {$\scriptstyle t_4$} [rb] <-1.3mm,0mm> at 0 2
		\put {$\scriptstyle t_5$} [rb] <-1.3mm,0mm> at 1 1
		\put {$\scriptstyle t_6$} [lb] <-4mm,0.4mm> at 2 0
		\put {$\scriptstyle t_7$} [lb] <-4mm,0mm> at 0 3
		\put {$\scriptstyle t_8$} [lb] <-4mm,0mm> at 1 2
		\put {$\scriptstyle b_1$} [lb] <0.8mm,0mm> at 2 1
		\put {$\scriptstyle b_2$} [lb] <0.8mm,0.8mm> at 3 0
	        \put {$\scriptstyle b_3$} [lb] <0.8mm,0.8mm> at 0 4
	        \put {$\scriptstyle b_4$} [lb] <0.8mm,0.8mm> at 1 3
	        \put {$\scriptstyle b_5$} [lb] <0.8mm,0.8mm> at  2 2
		\put {$\times$} at 0 0
	
		\put {$\circ$} at 3 0
		\put {$\circ$} at 2 2
		\put {$\circ$} at 2 1
		\put {$\circ$} at 1 3
		\put {$\circ$} at 0 4
		\endpicture}
\end{center}
Thus $\QQ[C] =\QQ[c_{11}, \dots, c_{85}]$ is a polynomial ring in~$40$ indeterminates. 
Notice that the dimension of~$\BO$ is $\dim(\BO) = \mu n = 16$, and that there are 
32 neighbour generators of the ideal $I(\BO)$. The linear parts of these generators are 
$$
\begin{array}{l}
c_{65},\;  c_{51} - c_{85},\;  c_{45},\;   c_{44},\;  c_{55},\;  c_{43} - c_{54},\;   
c_{42},\;  c_{41} - c_{75},\;  c_{52} - c_{75},\;  c_{35},  \cr 
c_{34},\;  c_{33},\;  c_{31},\;  c_{25},\;  c_{24},\;  c_{23},\;  c_{22},\;  c_{21},\;  
c_{32},\;  c_{15},\;  c_{14},\;  c_{13},\;  c_{12},\;   c_{11} .
\end{array}
$$
If we let $U$ be the union of the supports of these elements, we get 
$$
C\setminus U = \{c_{53},\,  c_{61},\,  c_{62},\,  c_{63},\,  c_{64},\,  c_{71},\,  
c_{72},\,  c_{73},\,  c_{74},\,  c_{81},\,  c_{82},\,  c_{83},\,  c_{84} \}
$$
which is exactly the set of basic indeterminates by Lemma~\ref{lem-equivclasses}.b.
Moreover, note that $C\setminus U$ is contained in the set 
of rim indeterminates (see Theorem~\ref{thm-RimIndets}.a).

For the trivial cotangent equivalence class~$E_0$ and the 
proper cotangent equivalence classes $E_1,\dots,E_q$, we get
\begin{align*}
E_0 = & \{ c_{11},  c_{12},  c_{13},  c_{14},  c_{15},  c_{21},  c_{22},  c_{23},  
c_{24},  c_{25},  c_{31},  c_{32},  c_{33},  c_{34},  c_{35}, \\
 &\; c_{42},  c_{44},  c_{45},  c_{55},  c_{65} \}\\
E_1 = & \{c_{51},  c_{85}\},\qquad  E_2 = \{ c_{43}, c_{54}\},\qquad
E_3 = \{c_{41}, c_{52}, c_{75}\} .
\end{align*}
Using $\sigma= {\tt DegRevLex}$, we obtain $E_1^\sigma = \{c_{51}\}$, $E_2^\sigma = \{c_{43}\}$,
and $E_3^\sigma = \{c_{41},  c_{52} \}$.
Next we compute the set $S_\sigma$ of the minimal generators of the leading term ideal 
of~$I_L = \langle \Lin_{\M}(I(\BO)) \rangle$ and get
$$
\begin{array}{l}
S_\sigma = \{  c_{11},  c_{12},  c_{13},  c_{14},  c_{15},  c_{21},  c_{22},  c_{23},  
c_{24},  c_{25},  c_{31},  c_{32},  c_{33},  c_{34},  c_{35}, \cr  
\qquad \quad c_{41}, c_{42},  c_{43},  c_{44},  c_{45},  c_{51},  c_{52},  c_{55},  c_{65}
  \} 
\end{array}
$$
which coincides with $E_0\cup E_1^\sigma \cup E_2^\sigma \cup E_3^\sigma$ 
(see the first claim of Theorem~\ref{thm-shapeofSsigma}.a).
Here we have $\#S_\sigma = 24$, $\#E_0 = 20$, $\#E_1 = 2$, $\#E_2 = 2$,  $\#E_3 = 3$,
and therefore 
$$
\#S_\sigma = \#E_0 + \#E_1 +\#E_2 +\#E_3 -3
$$
in accordance with the second claim of Theorem~\ref{thm-shapeofSsigma}.a.

The minimal sets of terms generating the sets in $\LTGFan(I_L)$ are
\begin{align*}
	Z_1 = & E_0\cup \{ c_{51} \} \cup \{ c_{43} \} \cup \{ c_{41},\ c_{52}\} \\
	Z_2 = & E_0\cup \{ c_{51} \} \cup \{ c_{43} \} \cup \{ c_{41},\ c_{75}\} \\
	Z_3 = & E_0\cup \{ c_{85} \} \cup \{ c_{43} \} \cup \{ c_{41},\ c_{52}\} \\
	Z_4 = & E_0\cup \{ c_{85} \} \cup \{ c_{43} \} \cup \{ c_{41},\ c_{75}\} \\
	Z_5 = & E_0\cup \{ c_{51} \} \cup \{ c_{54} \} \cup \{ c_{41},\ c_{52}\} \\
	Z_6 = & E_0\cup \{ c_{51} \} \cup \{ c_{54} \} \cup \{ c_{41},\ c_{75}\} \\
	Z_7 = & E_0\cup \{ c_{85} \} \cup \{ c_{54} \} \cup \{ c_{41},\ c_{52}\} \\
	Z_8 = & E_0\cup \{ c_{85} \} \cup \{ c_{54} \} \cup \{ c_{41},\ c_{75}\} \\
	Z_9 = & E_0\cup \{ c_{51} \} \cup \{ c_{43} \} \cup \{ c_{52},\ c_{75}\} \\
	Z_{10} = & E_0\cup \{ c_{85} \} \cup \{ c_{43} \} \cup \{ c_{52},\ c_{75}\} \\
	Z_{11} = & E_0\cup \{ c_{51} \} \cup \{ c_{54} \} \cup \{ c_{52},\ c_{75}\} \\
	Z_{12} = & E_0\cup  \{ c_{85} \} \cup \{ c_{54} \} \cup \{ c_{52},\ c_{75}\} .
\end{align*}
Thus we see that $\prod_{i=1}^3 \#E_i =12= \# \LTGFan(I_L)$
(see Theorem~\ref{thm-shapeofSsigma}.d).
As remarked before, $C\setminus U$ is the set of basic indeterminates, and we notice  
that $Z_i\cap (C\setminus U) =\emptyset$ for $i=1,\dots,12$ (see Theorem~\ref{thm-RimIndets}.a).

The complements of the sets~$Z_i$ in~$C$ are 
$$
\begin{array}{ll}
Y_1&\!\!\!= \{ c_{53}, c_{54}, c_{61}, c_{62}, c_{63}, c_{64}, c_{71}, c_{72}, 
c_{73}, c_{74}, c_{75}, c_{81}, c_{82}, c_{83}, c_{84}, c_{85} \} \cr
Y_2&\!\!\!= \{  c_{52}, c_{53}, c_{54}, c_{61}, c_{62},  c_{63},  c_{64},  
c_{71},  c_{72},  c_{73},  c_{74},  c_{81},  c_{82},  c_{83},  c_{84},  
c_{85}  \} \cr
Y_3&\!\!\!= \{ c_{51},  c_{53},  c_{54},  c_{61},  c_{62},  c_{63},  c_{64},  
c_{71},  c_{72},  c_{73},  c_{74},  c_{75},  c_{81},  c_{82},  c_{83},  
c_{84} \} \cr
Y_4&\!\!\!= \{ c_{51},  c_{52},  c_{53},  c_{54},  c_{61},  c_{62},  c_{63},  
c_{64},  c_{71},  c_{72},  c_{73},  c_{74},  c_{81},  c_{82},  c_{83},  
c_{84} \} \cr
Y_5&\!\!\!= \{ c_{43},  c_{53},  c_{61},  c_{62},  c_{63},  c_{64},  c_{71},  
c_{72},  c_{73},  c_{74},  c_{75},  c_{81},  c_{82},  c_{83},  c_{84},  
c_{85} \} \cr
Y_6&\!\!\!=\{ c_{43},  c_{52},  c_{53},  c_{61},  c_{62},  c_{63},  c_{64},  
c_{71},  c_{72},  c_{73},  c_{74},  c_{81},  c_{82},  c_{83},  c_{84},  
c_{85}  \} \cr
Y_7&\!\!\!=\{ c_{43},  c_{51},  c_{53},  c_{61},  c_{62},  c_{63},  c_{64},  
c_{71},  c_{72},  c_{73},  c_{74},  c_{75},  c_{81},  c_{82},  c_{83},  
c_{84} \} \cr
Y_8&\!\!\!= \{ c_{43},  c_{51},  c_{52},  c_{53},  c_{61},  c_{62},  c_{63},  
c_{64},  c_{71},  c_{72},  c_{73},  c_{74},  c_{81},  c_{82},  c_{83},  
c_{84}  \} \cr
Y_9&\!\!\!=\{ c_{41},  c_{53},  c_{54},  c_{61},  c_{62},  c_{63},  c_{64},  
c_{71},  c_{72},  c_{73},  c_{74},  c_{81},  c_{82},  c_{83},  c_{84},  
c_{85}  \} \cr
Y_{10}&\!\!\!=\{ c_{41},  c_{51},  c_{53},  c_{54},  c_{61},  c_{62},  c_{63},  
c_{64},  c_{71},  c_{72},  c_{73},  c_{74},  c_{81},  c_{82},  c_{83},  
c_{84}  \} \cr
Y_{11}&\!\!\!=\{ c_{41},  c_{43},  c_{53},  c_{61},  c_{62},  c_{63},  c_{64},  
c_{71},  c_{72},  c_{73},  c_{74},  c_{81},  c_{82},  c_{83},  c_{84},  
c_{85} \} \cr
Y_{12}&\!\!\!= \{ c_{41},  c_{43},  c_{51},  c_{53},  c_{61},  c_{62},  
c_{63},  c_{64},  c_{71},  c_{72},  c_{73},  c_{74},  c_{81},  c_{82},  
c_{83},  c_{84} \} .
\end{array}
$$
It is straightforward to verify that $C\setminus U$, the set of basic indeterminates, 
is contained in each set~$Y_i$ (see Theorem~\ref{thm-classifycij}.a).
Since we have $\#Z_i = \dim_\QQ(\Lin_\M(I(\BO)))$ for every~$i=1,\dots 12$, we are in the 
situation considered in Theorem~\ref{thm-classifycij}.c.
Using Algorithm~\ref{alg-compEmb}, we check that for each~$Z_i$ there exists 
an optimal $Z_i$-separating re-embedding of~$I(\BO)$.
Finally, for $i\in \{1,\dots,12\}$, we use $\#Z_i = 24$ and conclude 
that the set~$Z_i$ defines an isomorphism $\BO\cong \QQ[Y_i]$, 
where $\QQ[Y_i]$ is a polynomial ring having 16 indeterminates.
\end{example}

%%%%%%%%%%%%%%%%%%%%%%%%%%%%%%%%
%
%   Acknowledgements
%
%%%%%%%%%%%%%%%%%%%%%%%%%%%%%%%%%
\section*{Acknowledgements}

The third author thanks the University of Passau for its hospitality
during part of the preparation of this paper. We thank the referees
of the paper for their careful and detailed comments.

%%%%%%%%%%%%%%%%%%%%%%%%%%%%%%%%%%%%%%%%%%%%%%%%%%%%%%%
%
%   Bibliography
%
%%%%%%%%%%%%%%%%%%%%%%%%%%%%%%%%%%%%%%%%%%%%%%%%%%%%%%%

\end{document}